\documentclass[3p, 11pt]{elsarticle}

\makeatletter
\def\ps@pprintTitle{%
 \let\@oddhead\@empty
 \let\@evenhead\@empty
 \def\@oddfoot{}%
 \let\@evenfoot\@oddfoot}
\makeatother

\usepackage{lipsum}
\usepackage{hyperref}
\usepackage{amssymb}
\usepackage{caption,subcaption}
\usepackage{amsmath}
\usepackage{pgfplots}
\usepackage{pgf}
\usepackage{tikz} 
\usepackage{amsfonts}
\usepackage{graphicx}
\usepackage{epstopdf}
\usepackage{MnSymbol}
\usepackage{overpic}
\usepackage{float}
\usepackage{graphicx}
\usepackage[absolute,overlay]{textpos}
\usepackage{booktabs}
\usepackage{relsize}
\usepackage{mathrsfs}
\newtheorem{remark}{Remark}
\usepackage[linesnumbered, ruled, vlined]{algorithm2e}
\biboptions{sort&compress}
\newcommand{\bm}[1]{\ensuremath{\mathbf{#1}}}

\journal{Computer Methods in Applied Mechanics and Engineering}

\begin{document}

\begin{frontmatter}

\title{A DEIM Tucker Tensor Cross Algorithm  and its Application to Dynamical Low-Rank Approximation}

\author{Behzad Ghahremani, Hessam Babaee$^*$}

\address{Department of Mechanical Engineering and Materials Science, University of
Pittsburgh, 3700 O’Hara Street, Pittsburgh, PA, 15213, USA \\ \vspace{2mm}
* Corresponding Author, Email:h.babaee@pitt.edu \vspace{-8mm}}

\begin{abstract} 

We introduce a Tucker tensor cross approximation method that constructs a low-rank representation of a $d$-dimensional tensor by sparsely sampling its fibers. These fibers are selected using the discrete empirical interpolation method (DEIM). Our proposed algorithm is referred to as DEIM fiber sampling (\texttt{DEIM-FS}). For a rank-$r$ approximation of an $\mathcal{O}(N^d)$ tensor, \texttt{DEIM-FS} requires access to only $dNr^{d-1}$ tensor entries, a requirement that scales linearly with the tensor size along each mode. We demonstrate that \texttt{DEIM-FS}  achieves an approximation accuracy close to the Tucker-tensor approximation obtained via higher-order singular value decomposition at a significantly reduced cost. We also present \texttt{DEIM-FS (iterative)} that does not require access to singular vectors of the target tensor unfolding and can be viewed as a black-box Tucker tensor algorithm.  We employ \texttt{DEIM-FS} to reduce the computational cost associated with solving nonlinear tensor differential equations (TDEs) using dynamical low-rank approximation (DLRA). The computational cost of solving DLRA equations can become prohibitive when the exact rank of the right-hand side tensor is large. This issue arises in many TDEs, especially in cases involving non-polynomial nonlinearities, where the right-hand side tensor has full rank. This necessitates the storage and computation of tensors of size $\mathcal{O}(N^d)$. We show that \texttt{DEIM-FS}  results in significant computational savings for DLRA by constructing a low-rank Tucker approximation of the right-hand side tensor on the fly. Another advantage of using \texttt{DEIM-FS} is to significantly simplify the implementation of DLRA equations, irrespective of the type of TDEs. We demonstrate the efficiency of the algorithm through several examples including solving high-dimensional partial differential equations. 
\end{abstract}

\begin{highlights}
\item A novel cross Tucker tensor algorithm is introduced aimed at constructing near-optimal low-rank Tucker tensor models by sampling a relatively small number of elements from the target tensor.
\item The algorithm strategically samples fibers using the discrete empirical interpolation method.
\item An iterative cross algorithm is introduced, which operates without the need for accessing singular vectors of the tensor unfolding. It can be regarded as a black-box algorithm for constructing Tucker tensor models.
\item The cross algorithm is utilized to decrease both the computational cost and memory requirements when solving high-dimensional nonlinear tensor differential equations using dynamical low-rank approximation.   
\end{highlights}


\end{frontmatter}

\section{Introduction}
\label{sec:Introduction}
Multi-dimensional tensors play a critical role in many applications in science and engineering \cite{KB09}. However, performing computational tasks involving high-dimensional tensors or even storing them suffers from the curse of dimensionality: for a tensor of size $N_1\times N_2\times \dots \times N_d $, the number of its elements increases exponentially as $d$ grows, resulting in $\mathcal{O}(N^d)$ elements, where $N=N_i$, $i=1,\dots, d$ is assumed. Various tensor low-rank approximations have been developed to mitigate this issue by leveraging multi-dimensional correlations \cite{GKT13}. These dimension reduction techniques aim to decrease the total number of tensor elements while allowing a controllable loss of accuracy. Some of the most common tensor low-rank approximation schemes include Tucker tensor decomposition \cite{T66}, CANDECOMP/PARAFAC (CP) \cite{CP70}, hierarchical Tucker tensor decomposition \cite{G10}, and tensor train decomposition \cite{TT11}.

The Tucker tensor low-rank approximation reduces the total number of elements of a $d$-dimensional tensor from $N^d$ to $r^d + rdN$, where $r \ll N$ represents the rank of the unfolded tensors along each mode. The Tucker tensor low-rank approximation is the building block for many tensor low-rank approximations and it finds applications in various fields, including computer vision, deep neural networks, data mining, numerical analysis, neuroscience, and more. For an excellent review of these applications, refer to \cite{TKTensor}.

One of the applications of the Tucker tensor decomposition is to reduce the computational cost of solving multi-dimensional partial differential equations (PDEs) using dynamical low-rank approximation (DLRA) \cite{OKDynamical}. Discretizing these PDEs in all dimensions except time results in tensor differential equations (TDEs) in the form of $d\mathcal{V}/dt=\mathscr{F}(\mathcal{V})$, where $\mathcal{V} \in \mathbb{R}^{N_1\times N_2 \times \dots \times N_d}$ and $\mathscr{F}(\mathcal{V})$ is the right had side tensor of the same size as $\mathcal{V}$.  Example applications include the Schrodinger equation \cite{Beck:2000aa}, the Fokker-Planck equation \cite{HRFokker}, the Boltzmann transport equation \cite{ABTensor}, and the Hamilton-Jacobi-Bellman equations \cite{SDTensor}, among others.  Solving these TDEs, even in moderate dimensions ($3<d<7$), using traditional numerical methods such as finite difference and finite element methods, encounters the issue of the curse of dimensionality \cite{IGTensor}.  DLRA mitigates this issue by solving TDEs on the manifold of low-rank Tucker tensors, in which explicit evolution equations for a low-rank Tucker model are obtained.  DLRA for Tucker tensor form is the extension of DLRA for matrix differential equations \cite{KL07}, which has found many diverse applications such as kinetics  \cite{EL18,HW22}, linear sensitivity analyses \cite{DCB22} and species transport equations in turbulent combustion \cite{RNB21}. 

Reducing the computational cost of solving  DLRA equations is the primary motivating application for the developments presented in this work. The  DLRA evolution equations involve the $\mathscr F(\mathcal{V})$ tensor. When the exact rank of $\mathscr F(\mathcal{V})$ is high, the computational cost of DLRA increases.  This occurs for linear TDEs with a large number of right-hand side terms and for TDEs with high-order polynomial nonlinearities. In the case of TDEs with nonpolynomial nonlinearities, such as exponential or fractional nonlinearity, the computational cost of solving the DLRA evolution equations exceeds that of the full-order model (FOM). This increased cost arises in TDEs with nonpolynomial nonlinearity because $\mathscr{F}(\mathcal{V})$ is a full-rank tensor even when $\mathcal{V}$ is low rank. Explicit formation of this tensor necessitates memory and floating-point operation costs similar to those of the FOM, i.e., $\mathcal{O}(N^d)$. 

One solution to mitigate the issue of the cost of DLRA for nonlinear TDEs is to approximate $\mathscr{F}(\mathcal{V})$ with another Tucker tensor approximation. This idea was recently successfully employed for solving nonlinear matrix differential equations (MDEs) on low-rank matrix manifolds \cite{MNAdaptive,DPNFB23}. As demonstrated in Section \ref{sec:DEIM-FS DLRA}, once $\mathscr{F}(\mathcal{V})$ is approximated in the Tucker form, it can be efficiently incorporated into the DLRA evolution equations. 

To maintain the computational advantages of DLRA, any competitive Tucker tensor decomposition algorithm should satisfy stringent accuracy-versus-cost criteria: (i) the algorithm should be fast—preferably depending linearly on $N$ in terms of floating-point operations (flops) and memory requirements; (ii) the algorithm should be accurate, as the error introduced by the Tucker low-rank approximation of the right-hand side tensor should not exceed the DLRA errors, which are typically very small, often within the range of  $\mathcal{O}(10^{-6})$ to $\mathcal{O}(10^{-10})$ in relative errors, depending on the DLRA rank and the local temporal integration error. To this end,  we review various techniques for computing Tucker tensor approximation and evaluate them based on the aforementioned criteria.

 Determining the \emph{best} Tucker tensor low-rank approximation of a tensor lacks a known closed solution. The higher-order singular value decomposition (HOSVD) is a reliable approach for computing a near-optimal Tucker tensor decomposition  \cite{DDV00}. The computational cost of computing the Tucker tensor decomposition using HOSVD scales at least linearly with the total number of elements of a tensor, i.e., $\mathcal{O}(N^d)$.   For example, a six-dimensional probability density function with $N=100$ grid points in each dimension results in a tensor with one trillion ($10^{12}$) elements. This cost is prohibitive for many applications --- certainly for DLRA --- when $N^d$ is very large.  The computational expense stems from two primary sources: 
\begin{enumerate}
\item \textbf{Computing the SVD:} The HOSVD algorithm requires performing the SVD of large matrices to compute the factor matrices. Specifically, it involves computing the SVD of $d$ matrices of size $N \times N^{d-1}$ obtained by unfolding the tensor along its $d$ modes. The core is then computed via an orthogonal projection of the tensor onto factor matrices. The computational cost of performing $d$ SVD scales with $\mathcal{O}(dN^{(d+1)})$. 

\item \textbf{Data access:} The HOSVD requires access to all elements of the tensor for computing the factor matrices and the core tensor. When dealing with very large values of $N^d$, it might be impossible to hold the entire tensor in memory due to resource limitations. Consequently, multiple loading of the entire data in chunks from the disk becomes necessary, which is significantly slower compared to loading from memory. For DLRA, computing any tensor element requires applying $\mathscr{F}(\sim)$ to $\mathcal{V}$, implying that data access incurs flops costs in addition to the memory requirements.

\end{enumerate}
Numerous algorithms have been proposed to tackle these issues. For reducing the SVD cost, various approaches exist; for instance, randomized HOSVD \cite{SARandomized}, higher-order interpolatory decomposition \cite{ASHigher}, and sequentially truncated HOSVD \cite{VVM12}. However, all of these algorithms require access to all tensor entries. In contrast, cross algorithms tackle both of these issues simultaneously.

Cross algorithms are powerful techniques for constructing low-rank matrix and tensor approximations. The first cross approximation was proposed for matrix low-rank approximations \cite{GTZ97}. Cross algorithms are also known as \emph{pseudoskeleton} or \emph{CUR decomposition}. The simplest cross algorithm is an interpolatory CUR matrix decomposition, wherein a rank-$r$ matrix is constructed by sampling only $r$ columns and rows of the matrix.   The clear advantage of the interpolatory CUR algorithm is that the low-rank approximation does not require all the entries of the matrix. 
  However, the accuracy of the low-rank approximation obtained via CUR critically depends on the selection of rows and columns. As shown in \cite{GTZ97}, the accuracy of the CUR approximation depends on the determinant of the intersection submatrix, which is referred to as \emph{matrix volume}.  In particular, the columns and rows should be selected such that the matrix volume is maximized. Since this selection problem is NP-hard, several heuristic algorithms have been proposed including Maxvol \cite{SGHowto, SGTheMaximal}, Cross2D  \cite{DSPolilinear, ETIncomplete}, leverage score \cite{MMRandomized}, and discrete empirical interpolation method (DEIM) \cite{DSaDEIM} algorithms.

The cross algorithms have been extended to various tensor low-rank approximations, including Tucker tensor decomposition \cite{CCGeneralizing}, tensor train decomposition \cite{ISTtcross}, and hierarchical Tucker tensor decomposition \cite{JBLinear}. Similar to the interpolatory CUR algorithm, these algorithms do not require access to all entries of the tensor. We refer the readers to \cite{SACross} for an excellent review of various tensor cross approximations.

The cross approximation, introduced in \cite{CCGeneralizing}, is known as fiber sampling Tucker decomposition (FSTD). FSTD is an elegant extension of matrix CUR decomposition to Tucker tensor approximations.   However,  as we demonstrate, FSTD becomes ill-conditioned as rank increases and therefore, the approximation error cannot be reduced to machine precision. The issue of ill-conditioning is similarly encountered in matrix CUR decompositions, where the intersection submatrix becomes singular as rank increases. This issue is mitigated by applying the QR factorization to either the selected columns or rows \cite{ISTtcross}.  As a result, FSTD is not well-suited for DLRA due to the stringent accuracy-versus-cost requirements that DLRA demands.

In this paper, we present a novel Tucker cross algorithm that addresses the aforementioned challenges.  In particular, the contributions of this paper are the following: 
\begin{enumerate}
    \item  We introduce \texttt{DEIM-FS} --- a Tucker tensor cross algorithm by sampling $r^{d-1}$ fibers along each mode of a tensor. Therefore, it requires access to $dr^{d-1}N$ entries of the tensor. The fiber selection (FS) is guided by the DEIM algorithm. The DEIM algorithm requires access to the exact or approximate singular vectors of the tensor unfolding along each mode. As we will discuss in this paper, DLRA is one such application, where the approximate singular vectors are available at no additional cost. We demonstrate that the low-rank approximation error of the \texttt{DEIM-FS} algorithm is comparable to that achieved by HOSVD. 
    \item  We present  \texttt{DEIM-FS~(iterative)} for the problems where the exact or approximate singular vectors are not available. \texttt{DEIM-FS (iterative)}  starts with a random guess for fibers and iteratively applies \texttt{DEIM-FS} until convergence. In practice, a small number of iterations are needed.  As a result, \texttt{DEIM-FS~(iterative)} can be regarded as a \emph{black-box} Tucker cross algorithm.  
    \item We present \texttt{DLRA-DEIM-FS}  to reduce the computational cost of solving DLRA evolution equations for nonlinear TDEs. \texttt{DLRA-DEIM-FS}  constructs a low-rank Tucker tensor approximation of the right-hand side of the TDE by applying \texttt{DEIM-FS} to  $\mathscr{F}(\mathcal{V})$.  We augment the  \texttt{DEIM-FS} with \emph{rank-adaptivity}, where the Tucker tensor rank is adjusted on the fly to meet an error threshold criterion.  
\end{enumerate}

The paper is organized as follows. The methodologies are discussed in
Section \ref{sec:Method}, demonstrations and results are presented in Section \ref{sec:Demo}, and the conclusions follow in Section \ref{sec:conclusions}.

\section{Methodology}
\label{sec:Method}

\subsection{Definitions and notations} \label{sec:Definition}
We first introduce the notation used for vectors, matrices, and tensors. Vectors are denoted in bold lowercase letters (e.g. $\bm a$),  matrices are denoted by bold uppercase letters (e.g. $\bm A$), and tensors by uppercase calligraphic letters (e.g. $ \mathcal{F}$). The symbol $\bigtimes_n$ is used to denote the $n$-mode product. The $n$-mode product of a tensor $\mathcal{F} \in \mathbb{R}^{N_1 \times N_2 \times \dots \times N_d}$ with a matrix $\bm M \in \mathbb{R}^{J \times N_n}$ is obtained by $\mathcal{F} \bigtimes_n \bm M$ and is of size $N_1 \times \dots \times N_{n-1} \times J \times N_{n+1} \times \dots \times N_d$. We denote the unfolding of tensor $\mathcal F$ along its $n$-th mode with $\mathcal{F}_{(n)}$. Unfolding a tensor involves reshaping its elements in such a way that results in a matrix instead of a tensor. For instance, a tensor of size $3 \times 5 \times 6$ can be unfolded along the second axis as a matrix of size $5 \times 18$ \cite{TKTensor}. The Frobenius norm of a tensor is shown by $\| \mathcal{F} \|_F$ and is defined as:
\begin{equation}\label{eq:FrobNorm}
    \| \mathcal{F} \|_F =  \sqrt{\sum_{n_1=1}^{N_1} \sum_{n_2=1}^{N_2} \dots \sum_{n_d=1}^{N_d} a^{2}_{n_1 n_2 \dots n_d}},
\end{equation}
where $ a_{n_1 n_2 \dots n_d}$ are the entries of tensor $\mathcal{F}$. We use typewriter font to denote algorithms, e.g., \texttt{SVD} or \texttt{DEIM-FS}.  
We use the MATLAB indexing notation where $\bm A(\bm p,:)$ selects all columns at the $\bm p$ rows and $\bm A(:, \bm p)$ selects all rows at the $\bm p$ columns of matrix $\bm A$ and $\bm p=[p_1, p_2, \dots, p_q]$ is the integer vector containing the selected indices. We also use MATLAB notation for computing the SVD of a matrix. For example, consider for $\bm A \in \mathbb{R}^{m\times n}$. Then $[\bm U,\bm \Sigma, \bm 
 V] =\texttt{SVD}(\bm A,r)$ means computing the SVD of $\bm A$ and truncating at rank $r$, where $r<m$ and $r<n$, $\bm U \in \mathbb{R}^{n\times r}$ is the matrix of left singular vectors,  $\bm \Sigma \in \mathbb{R}^{r \times r}$ is the matrix of singular values and $\bm V \in \mathbb{R}^{m \times r}$ is the matrix of right singular vectors.  If any of the singular matrices are not needed, the symbol ($\sim$) is used. For example, $[\bm U,\sim, \sim] =\texttt{SVD}(\bm A,r)$ returns only the first $r$ left singular vectors.

\subsection{Dynamical low-rank approximation for Tucker tensors} \label{sec:LowRankApprox}
As mentioned in the Introduction, DLRA is the primary motivation for the developments in this paper. The DLRA formulation for Tucker tensors is used to solve high-dimensional PDEs on the manifold of low-rank Tucker tensors \cite{OKDynamical}. In the following, we briefly introduce DLRA and explain the computational cost issues related to DLRA. 

We consider a general PDE given by:
\begin{equation}\label{eq:ContinousPDE}
    \frac{\partial v (\bm x,t)}{\partial t} = f\bigl( v(\bm x,t) \bigr),
\end{equation}
augmented with appropriate initial and boundary conditions. Here $\bm x = (x_1, x_2, \dots , x_d) \in \mathbb{R}^d$,~ $t$ is time, $f$ is a general nonlinear differential operator, and $d$ is the dimension of the problem. We consider the differential operators in $\bm x$ being discretized using a method of line. We discretize the differential operators of  Eq. (\ref{eq:ContinousPDE}) in $\bm x$ using a method of lines, which  results in the following TDE:
\begin{equation}\label{eq:TDE}
    \frac{d \mathcal{V}}{dt} = \mathscr{F}\bigl( \mathcal{V} \bigr),
\end{equation}
where $\mathcal{V}(t) \in \mathbb{R}^{N_1 \times N_2 \times \dots \times N_d}$ is the solution tensor and $\mathscr{F}$ is the discrete representation of $f$. Here, $N_1, N_2, \dots, N_d$ are the number of the discretized points along each mode of the tensor. We refer to Eq. \ref{eq:TDE}  as the FOM.

Discretizing Eq.~(\ref{eq:ContinousPDE}) via classical methods such as finite-difference and spectral methods results in a system where the degrees of freedom grow exponentially fast as the number of dimensions grows.  One approach to mitigate this issue is to solve Eq.~(\ref{eq:TDE}) on a manifold of low-rank Tucker tensors. To this end, consider a low-rank Tucker tensor approximation of $\mathcal{V}$  \cite{OKDynamical} as shown below:
\begin{equation}\label{eq:low-rankHOSVD}
    \mathcal{V}(t) \approx \hat{\mathcal{V}}(t) = \mathcal{S}(t) \times_1 \bm U^{(1)}(t) \times_2 \bm U^{(2)}(t) \dots \times_d \bm U^{(d)}(t) ,
\end{equation}
where $\times_n$ is tensor mode product, $\mathcal{S}\in \mathbb{R}^{r_1 \times r_2 \times \dots \times r_d}$ is the core tensor, $\bm U^{(i)} \in \mathbb{R}^{N_i \times r_i}$  are the orthonormal time-dependent bases or factor matrices along the corresponding mode of the tensor, i.e., ${\bm U^{(i)}}^T \bm U^{(i)} = \bm I$, where $\bm I$ is the identity matrix and  $r_i \ll N_i,  i = 1,2, \dots, d$ are the rank along each tensor mode. Substituting the low-rank approximation given by Eq.~(\ref{eq:low-rankHOSVD})   into Eq.~(\ref{eq:TDE}) results in a residual equal to:
\begin{multline}\label{eq:minimize}
    \mathcal{R} (\Dot{\mathcal{S}}, \Dot{\bm U}^{(1)}, \Dot{\bm U}^{(2)}, \dots ,\Dot{\bm U}^{(d)}) = \\ 
  \Biggl \| \frac{d~\bigl(\mathcal{S} \times_1 \bm U^{(1)} \times_2 \bm U^{(2)}~ \dots \times_d \bm U^{(d)} \bigr)}{dt} - \mathscr{F} \bigl(\mathcal{S} \times_1 \bm U^{(1)} \times_2 \bm U^{(2)}~ \dots \times_d \bm U^{(d)} \bigr) \Bigg \|_F^2.
\end{multline}
The evolution equations for orthonormal bases and the core tensor are  obtained by minimizing  the above residual 
subject to the orthonormality constraints of the bases, which results in the following evolution equations for the core and the factor matrices \cite{OKDynamical}: 
\begin{subequations}
\begin{align}
  \Dot{\mathcal{S}}  &= \mathcal{F} ~ \bigtimes_{i=1}^d ~ {\bm U^{(i)}}^T ,
  \label{eq:Sevolution} \\
  \Dot{\bm U}^{(i)}  &= \Bigl(\bm I - \bm U^{(i)} \bm U^{(i)^T} \Bigr) ~ \Bigl[ \mathcal{F} \bigtimes_{k \neq i} {\bm U^{(k)}}^T \Bigr]_{(i)} ~ \mathcal{S}_{(i)}^\dag ,
  \label{eq:Uevolution}
\end{align}
\end{subequations}
where $\bm I$ is the identity matrix, $\mathcal{S}_{(i)}^\dag = \mathcal{S}_{(i)}^T \bigl(\mathcal{S}_{(i)}\mathcal{S}_{(i)}^T)^{-1}$ is the pseudo-inverse, $\mathcal{F} \in \mathbb{R}^{N_1 \times N_2 \times \dots \times N_d}$ is a tensor defined as $\mathcal{F} = \mathscr{F} \bigl(\mathcal{S} \times_1 \bm U^{(1)} \times_2 \bm U^{(2)}~ \dots \times_d \bm U^{(d)}\bigr)$. 
For further details on the derivation of the evolution equations,  refer to \cite{OKDynamical}. Eqs. (\ref{eq:Sevolution}) and (\ref{eq:Uevolution}) represent the DLRA evolution equations in the Tucker tensor form. For recent developments related to a stable time integration scheme and rank adaptivity of DLRA equations,  see \cite{CL20,CKL22}.

The computational advantage of the DLRA evolution equations over the FOM (Eq.~\ref{eq:TDE}) is that in Eqs.~(\ref{eq:Sevolution}) and (\ref{eq:Uevolution}),  the solution is sought in terms of the factor matrices  ($\bm U^{(i)}$) and the core tensor  ($\mathcal{S}$) instead of the full-dimensional tensor ($\mathcal{V}(t)$). The memory requirement for storing $\{\mathcal{S}, \bm U^{(i)}\}$ is $\mathcal{O}(r^d) + \mathcal{O}(rdN)$. However, solving $\mathcal{V}(t)$ using the FOM requires $\mathcal{O}(N^d)$ memory. For simplicity in computational complexity analysis, we assume \( r = r_1 = r_2 = \ldots = r_d \) and \( N = N_1 = N_2 = \ldots = N_d \). 

The computational savings of the DLRA in memory and floating-point operations (flops) are lost when dealing with TDEs featuring general nonlinearity. In these cases, the right-hand side (RHS) tensor \( \mathcal{F} \) is full rank, necessitating its computation and storage in memory. The computational cost of computing \( \mathcal{F} \) scales at \( \mathcal{O}(N^d) \), mirroring the computational complexity of solving the FOM. Even in linear TDEs, a substantial number of terms on the right-hand side may result in a large exact rank for \( \mathcal{F} \). Therefore, the issue of high computational cost is not only limited to general nonlinearity, but the issue arises for problems in which the exact rank of \( \mathcal{F} \) is large.  Consequently, evaluating \( \mathcal{F} \) remains the primary computational bottleneck for DLRA. To address this challenge, we introduce an efficient and innovative cross algorithm that constructs a low-rank Tucker tensor approximation of \( \mathcal{F} \) by selectively sampling a few fibers of \( \mathcal{F} \) using the DEIM algorithm.

In the next section, we first provide the utility of the DEIM algorithm for nonlinear reduced-order modeling. The cross algorithm presented in this paper is inspired by our previous work \cite{MNAdaptive}, where we developed a CUR algorithm for DLRA of nonlinear MDEs. In Section \ref{sec:DEIM-MDE}, we provide a brief overview of the algorithm presented in \cite{MNAdaptive}.

\subsection{DEIM for low-rank approximation of vector differential equations}  \label{sec:DEIM-Vec}
In reduced-order modeling based on proper orthogonal decomposition (POD), analogous challenges arise due to general nonlinearity. Consider the FOM given by: $d\bm{v}/dt=f(\bm{v})$, where $\bm{v} \in \mathbb{R}^N$ represents the state vector, and $f(\bm{v}) : \mathbb{R}^N \rightarrow \mathbb{R}^N$ is the RHS vector. Let $\bm{U} \in \mathbb{R}^{N\times r}$ denote the matrix of orthonormal POD modes, and $\bm{y} \in \mathbb{R}^r$ be the POD coefficients, such that $\hat{\bm{v}} = \bm{U} \bm{y} \in \mathbb{R}^N$ is the POD approximation of $\bm{v}$. Here, $r \ll N$ represents the number of POD modes. The ROM is obtained via Galerkin projection: $d\bm{y}/dt = \bm{U}^T f(\bm{U} \bm{y})$. If $f(\sim)$ involves a polynomial nonlinearity of degree $p$, it becomes feasible to compute $\bm{U}^T f(\bm{U} \bm{y})$ with a computational cost of at least $\mathcal{O}(r^p)$. This implies the potential avoidance of forming the vector $f(\bm{U} \bm{y}) \in \mathbb{R}^N$, thus preventing the computational cost of solving the FOM from scaling with the FOM size ($N$) \cite{MBKB18}. However, in scenarios where $f(\sim)$ has high-order polynomial nonlinearity (i.e., a large $p$), the computational cost of solving the ROM can become considerable. Furthermore, in cases where $f(\sim)$ has non-polynomial nonlinearity, the explicit formation of the vector $f(\bm{U} \bm{y})$ becomes necessary. This results in the loss of computational savings offered by the POD-ROM, as computing $f(\bm{U} \bm{y})$ requires $\mathcal{O}(N)$ operations—equivalent to solving the FOM.

One computationally efficient remedy is to interpolate the  $f(\bm U \bm y )$ onto a low-rank basis using the DEIM algorithm \cite{SCNonlinear}. This involves sampling $f(\bm U \bm y )$  at only a few strategically selected points.  To explain this algorithm, let $\bm f=f(\bm U \bm y ) $ and $\bm U_f \in \mathbb{R}^{N\times r_f}$ be a low-rank basis for the vector $\bm f$, where $r_f$ is the number of POD modes for the vector $\bm f$. The basis  $\bm U_f$ is computed in the offline stage as the left singular vectors of the RHS snapshot matrix. The DEIM algorithm \cite[Algorithm 1]{SCNonlinear} yields a set of near-optimal sampling points for interpolation of vector $\bm f$ onto $\bm U_f$:
\begin{equation}
    \bm p = \texttt{DEIM}(\bm U_f),
\end{equation}
where $\bm p =[p_1, p_2, \dots, p_r]$ is the vector of sampling point indices. For convenience, the DEIM algorithm is provided in Appendix 1. The vector $\bm f$ can be interpolated onto $\bm U_f$ using:
\begin{equation}\label{eq:DEIM-f}
    \hat{\bm f} = \bm{U}_f \bm{U}_f(\bm p, :)^{-1}\bm f(\bm p).
\end{equation}
It is easy to verify that  $\hat{\bm f}$ and $\bm f$ are equal to each other at interpolation points, i.e.,  $\hat{\bm f}(\bm p) = \bm f(\bm p)$. Incorporating Eq.~\ref{eq:DEIM-f} into the POD-ROM results in: $d\bm y/dt = \bm{U}^T\bm{U}_f\bm{U}_f(\bm p, :)^{-1}\bm f(\bm p)$. The small matrix $\bm{U}^T\bm{U}_f\bm{U}_f(\bm p, :)^{-1} \in \mathbb{R}^{r \times r_f}$ can be computed and stored in the offline stage.   The key advantage of using the DEIM algorithm in POD-ROM is that it requires evaluating $f(\bm U \bm y)$ at a only small number of points ($r_f$)- regardless of the type of nonlinearity of $f(\sim)$. As shown in \cite{SCNonlinear}, the approximation error of $\hat{\bm f}$ is bounded by the best approximation error by a magnification factor:
\begin{equation}
\|\hat{\bm f} - \bm f\| \leq \eta \|(\bm I - \bm U_f^T \bm U_f)\bm f \|
\end{equation}
where $\eta= \| \bm{U}_f(\bm p, :)^{-1} \|$ is the magnification factor and $\|(\bm I - \bm U_f^T \bm U_f)\bm f \|$ is the optimal error of approximating $\bm f$ in the span of $\bm U_f$, which is obtained via the orthogonal projection of $\bm f$ onto $\bm U_f$. The DEIM algorithm is designed to minimize $\eta$ using a greedy approach.

\subsection{DEIM  for low-rank approximation of matrix differential equations}  \label{sec:DEIM-MDE}
The DEIM algorithm, developed for the low-rank approximation of MDEs in \cite{MNAdaptive}, deals with an MDE expressed as $d\bm{V}/dt = F(\bm{V})$, where $\bm{V} \in \mathbb{R}^{N_1 \times N_2}$ and $F(\bm{V}): \mathbb{R}^{N_1 \times N_2} \rightarrow \mathbb{R}^{N_1 \times N_2}$. The DLRA aims to find solutions for the above MDE constrained to the manifold of rank $r$ matrices, where $r \ll N_1$ and $r \ll N_2$. DLRA for MDE can be derived as a specific case of DLRA for TDEs (refer to Eq.~\ref{eq:Sevolution}-\ref{eq:Uevolution}) by setting $d=2$ and $r=r_1=r_2$. When the exact rank of $F(\bm{V})$ is large, the computational cost of DLRA increases. For example,  when  $F(\sim)$ exhibits non-polynomial nonlinearity, $F(\bm{V})$ is full rank, and the explicit formation of $F(\bm{V})$ is needed. This results in the loss of computational savings provided by DLRA.

The algorithm proposed in \cite[Algorithm 1]{MNAdaptive} constructs a low-rank approximation of the matrix $\bm{F} = F(\bm{V})$ by selectively sampling $r$ columns and $r$ rows from $\bm{F}$. This involves:
\begin{enumerate}
\item  Sampling the columns of the RHS matrix: $\bm{F}(:, \bm{p}_2) \in \mathbb{R}^{N_1 \times r}$, where $\bm{p}_2 \in \mathbb{I}^r$ denotes the integer vector containing  the column indices.
\item Conducting the QR decomposition of matrix $\bm{F}(:, \bm{p}_2) = \bm{Q} \bm{R}$ to construct an orthonormal basis $\bm{Q}$ for the columns of matrix $\bm{F}$.
\item Computing the rows of the RHS matrix $\bm{F}(\bm{p}_1, :) \in \mathbb{R}^{N_1 \times r}$, where $\bm{p}_1 \in \mathbb{I}^r$ represents the integer row indices.
\item Interpolating each column of matrix $\bm{F}$ onto the basis $\bm{Q}$ using the computed values at rows indexed by $\bm{p}_1$: $\bm{F} \approx \hat{\bm{F}} = \bm{Q} \bm{Q}(\bm{p}_1,:)^{-1} \bm{F}(\bm{p}_1, :)$.
\end{enumerate}

In the above algorithm, the row and column indices $\bm p_1$ and $\bm p_2$ are determined using the DEIM algorithm: $\bm p_1=\texttt{DEIM}(\bm U^{(1)}_F)$ and $\bm p_2=\texttt{DEIM}(\bm U^{(2)}_F)$, where $\bm U^{(1)}_F$ and $\bm U^{(2)}_F$ are the left and right singular vectors of matrix $\hat{\bm F}$.

\subsection{DEIM  for  Tucker tensor cross approximation}  \label{sec:DEIM-FS algorithm}
Cross tensor approximations are extensions of matrix CUR approximation techniques, offering a practical approach for efficiently estimating low-rank tensors \cite{SACross}. In this paper, a novel cross tensor approximation technique is proposed. The main steps of the proposed methodology are summarized in Algorithm \ref{Algorithm1}. We refer to our algorithm as DEIM fiber sampling (\texttt{DEIM-FS}). This algorithm is presented for a three-dimensional tensor for simplicity, but it can easily be generalized to higher-dimensional tensors.  The algorithm constructs a low-rank Tucker tensor approximation of  $\mathcal{F}$ by only sampling a few fibers along each mode of the tensor.

To this end, consider a Tucker low-rank approximation of $\mathcal{F}$  given by:
\begin{equation}\label{eq:DEIM-FS-approx}
    \mathcal{F} \approx  \mathcal{S}_{\mathcal{F}} \times_1 {\bm U^{(1)}_{\mathcal{F}}} \times_2 {\bm U^{(2)}_{\mathcal{F}}} \times_3 {\bm U^{(3)}_{\mathcal{F}}},
\end{equation}
where $\mathcal{F} \in \mathbb{R}^{N_1 \times N_2 \times N_3}$, $\mathcal{S}_{\mathcal{F}} \in \mathbb{R}^{{r_{\mathcal{F}}}^{}_1 \times {r_{\mathcal{F}}}^{}_2 \times {r_{\mathcal{F}}}^{}_3}$, ${\bm U_{\mathcal{F}}}^{(1)} \in \mathbb{R}^{N_1 \times {r_{\mathcal{F}}}^{}_1}$, ${\bm U^{(2)}_{\mathcal{F}}} \in \mathbb{R}^{N_2 \times {r_{\mathcal{F}}}^{}_2}$, ${\bm U^{(3)}_{\mathcal{F}}} \in \mathbb{R}^{N_3 \times {r_{\mathcal{F}}}^{}_3}$, and $\bm r_{\mathcal {F}} = (r_{\mathcal {F}_1},r_{\mathcal {F}_2},r_{\mathcal {F}_3})$ is the multi-rank of the Tucker tensor decomposition.  The number of selected fibers along the first, second, and third modes are denoted by $ r'_{\mathcal {F}_1}, r'_{\mathcal {F}_2}, r'_{\mathcal {F}_3}$, respectively. As we will explain below (See Remark \ref{rm:rank}), the number of selected fibers must be greater than or equal to the target Tucker low-rank, i.e., $r'_{\mathcal {F}_i} \geq r_{\mathcal {F}_i}$. According to the numerical results of Section \ref{sec:ToyEx}, we demonstrate that $r'_{\mathcal {F}_i}=r_{{\mathcal{F}}_i} + 2$ is a reasonable choice and this choice is used in all demonstrations in this paper.  The indices of the selected fibers are shown by vectors $\bm p_1, \bm p_2, \bm p_3$, where $\bm p_i$ is an integer vector containing  $r'_{\mathcal {F}_i}$.  

\begin{figure}[t]
  \centering
  \label{DEIM-FS}\includegraphics[scale=0.45]{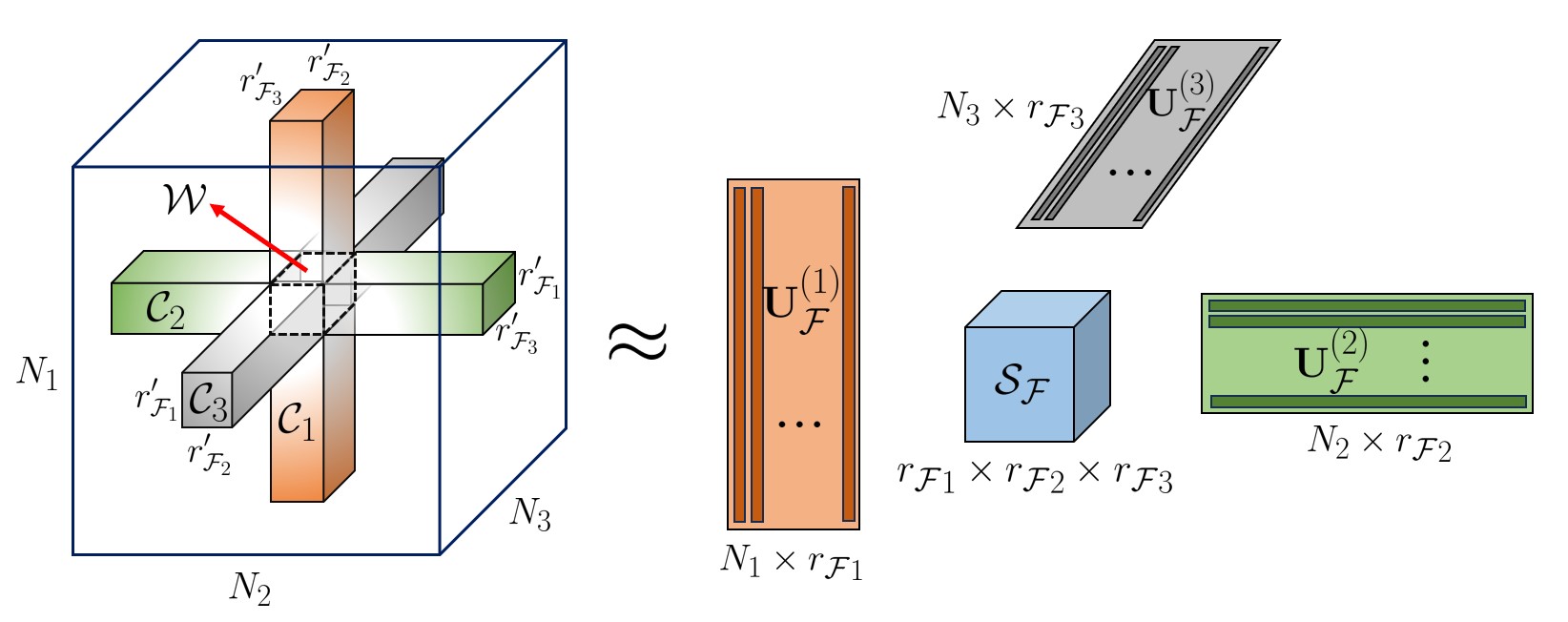}
  \caption{Schematic of the \texttt{DEIM-FS} cross algorithm for a 3D tensor. For simplicity, we assume that all selected fibers are adjacent to each other. }
  \label{fig:DEIM-FS}
\end{figure}

The schematic of the algorithm is shown in Figure \ref{fig:DEIM-FS}, where  $\mathcal{C}_1=\mathcal{F}(:, \bm p_2, \bm p_3) \in \mathbb{R}^{N_1 \times r'_{\mathcal {F}_2} \times r'_{\mathcal {F}_3}},~\mathcal{C}_2=\mathcal{F}( \bm p_1,:, \bm p_3) \in \mathbb{R}^{N_2 \times r'_{\mathcal {F}_1} \times r'_{\mathcal {F}_3}},~\mathcal{C}_3= \mathcal{F}(\bm p_1, \bm p_2, :) \in \mathbb{R}^{N_3 \times r'_{\mathcal {F}_1} \times r'_{\mathcal {F}_2}}$  are the sub-tensors formed by clustering the selected fibers of $\mathcal{F}$ along each direction.   For \texttt{DEIM-FS}, we assume that the first $r'_{\mathcal {F}_i}$ left singular vectors of matrix $\mathcal{F}_{(i)}$  or some close approximation of these vectors are known. We denote these singular vectors with $\tilde{\bm U}_{\mathcal{F}}^{(i)}\in \mathbb{R}^{N_i \times r'_{\mathcal {F}_i}}$, where   $i=1,2,3$ denotes different unfolding. In Section \ref{sec:iterative_DEIM-FS}, we present \texttt{DEIM-FS (iterative)}, where the matrices of left singular vectors are not needed.  In the following, we explain all the steps of  \texttt{DEIM-FS}.

\subsubsection{Computing the factor matrices}
The factor matrices are computed using the DEIM-selected fibers of tensor $\mathcal{F}$.
  To this end,  the DEIM algorithm is applied to the left singular vectors of $\mathcal{F}_{(i)}$ matrices:
\begin{equation*}
   \bm p_i= \texttt{DEIM} (\tilde{\bm U}_{\mathcal{F}}^{(i)}), \quad  i=1,2,3.
\end{equation*}
This generates the fiber indices for each mode of tensor $\mathcal{F}$.

  The DEIM-selected fibers are extracted from  tensor $\mathcal F$ and the resulting subtensors are unfolded  along  the corresponding modes to obtain the following matrices:
\begin{equation*}
   \bm C_1= \Bigl(\mathcal{F}(:, \bm p_2, \bm p_3)\Bigr)_{(1)},   \quad \bm C_2 = \Bigl(\mathcal{F}(\bm p_1, :, \bm p_3)\Bigr)_{(2)}, \quad \mbox{and} \quad
    \bm C_3 = \Bigl(\mathcal{F}(\bm p_1, \bm p_2, :)\Bigr)_{(3)},
\end{equation*}
where $\bm C_1 \in \mathbb{R}^{N_1 \times r'_{\mathcal {F}_2} r'_{\mathcal {F}_3}}$, $\bm C_2 \in \mathbb{R}^{N_2 \times r'_{\mathcal {F}_1} r'_{\mathcal {F}_3}}$, and $\bm C_3 \in \mathbb{R}^{N_3 \times r'_{\mathcal {F}_1} r'_{\mathcal {F}_2}}$.  In the next step, the factor matrices are computed along each mode of the Tucker tensor decomposition. To this end,  we perform the SVD of the $\bm C_i$ matrices and truncate at rank $r_{\mathcal{F}_i}$:
\begin{equation}\label{Eq:SVD-Fiber}
[\bm U^{(i)}_{\mathcal{F}}, \sim , \sim ] = \texttt{SVD}( \bm C_i, r_{\mathcal{F}_i}), \quad i=1,2,3,
\end{equation}
where $\bm U^{(i)}_{\mathcal{F}} \in \mathbb{R}^{N_i \times r_{\mathcal{F}_i}}$ is the matrix of left singular vectors.
The rank $r_{\mathcal{F}_i}$ can be either fixed a priori or determined adaptively based on accuracy requirements as explained in Section \ref{sec:rankadaptive}. The $\bm U^{(i)}_{\mathcal{F}}$ matrices computed in this step are the factor matrices of the cross Tucker tensor decomposition. 
\subsubsection{Computing the core tensor}
Since the factor matrices are already computed,  the optimal core tensor may be obtained via the \emph{orthogonal projection} of $\mathcal{F}$ onto the factor matrices. However, the orthogonal projection would require all entries of $\mathcal{F}$, which is undesirable.   The intersection tensor is denoted with  $\mathcal{W}$ (see Figure \ref{fig:DEIM-FS}) and it is equal to:
\begin{equation*}
\mathcal{W} = \mathcal{F}(\bm p_1, \bm p_2, \bm p_3)\in\mathbb{R}^{r'_{\mathcal {F}_1}\times r'_{\mathcal {F}_2}\times r'_{\mathcal {F}_3}}.
\end{equation*}
Note that the entries of the intersection tensor are a subset of the already-computed fibers in the previous step.  

The core tensor is calculated such that the difference between the cross Tucker tensor approximation and the actual values of the tensor is minimized at the DEIM intersection tensor ($\mathcal{W}$).  More specifically, we seek to find a core tensor such that 
\begin{equation*}
    E_{\mathcal{W}}(\mathcal{S}_{\mathcal{F}}) = \big \|\mathcal{W} - \mathcal{S}_{\mathcal{F}} \times_1 ~ \bm U^{(1)}_{\mathcal{F}}(\bm p_1, :) \times_2 ~ \bm U^{(2)}_{\mathcal{F}}(\bm p_2, :) \times_3 ~ \bm U^{(3)}_{\mathcal{F}}(\bm p_3, :) \big \|^2_F,
\end{equation*}
is minimized. The above minimization problem has a closed-form solution and it is computed via the least squares solution as shown below:
\begin{equation}\label{eq:core_lst}
    \mathcal{S}_{\mathcal{F}} = \mathcal{W} \times_1 ~ \bm U^{(1)}_{\mathcal{F}}(\bm p_1, :)^\dagger \times_2 ~ \bm U^{(2)}_{\mathcal{F}}(\bm p_2, :)^\dagger \times_3 ~ \bm U^{(3)}_{\mathcal{F}}(\bm p_3, :)^\dagger.
\end{equation}
The above procedure can be also viewed as an \emph{oblique projection} of $\mathcal{F}$ onto the factor matrices similar to matrix CUR decompositions \cite{DSaDEIM,DPNFB23}. 

We denote the resulting Tucker tensor decomposition obtained via  Algorithm \ref{Algorithm1} by:
\begin{equation}\label{eq:DEIM-FS-T}
    \hat{\mathcal{F}}= \mathcal{S}_{\mathcal{F}} \times_1 {\bm U^{(1)}_{\mathcal{F}}} \times_2 {\bm U^{(2)}_{\mathcal{F}}} \times_3 {\bm U^{(3)}_{\mathcal{F}}}.
\end{equation}
If rank adaptivity is not desired, no further computation is required and the output of the algorithm is the Tucker tensor decomposition with the user-specified multi-rank of $\bm r_{\mathcal{F}} = (r_{\mathcal{F}_1},r_{\mathcal{F}_2},r_{\mathcal{F}_3})$.
\subsubsection{Rank adaptivity}\label{sec:rankadaptive}
In certain applications, there could be a specific requirement for a low-rank approximation error. Consequently, adjusting the rank becomes necessary to meet a user-specified criterion. To address this requirement, we demonstrate that the cross tensor approximation algorithm can be made rank-adaptive with minor modifications.  The error criterion is application-dependent. We consider $\epsilon_i$ 
\begin{equation}\label{eq:ErrorProxy}
\epsilon_i = \dfrac{\min ({\bm \Sigma_{\mathcal{F}}}_i)}{\|{\bm \Sigma_{\mathcal{F}}}_i \|_F}  \quad i=1,2,3~,
\end{equation}
as the error proxy, where $\bm \Sigma_{\mathcal{F}_i}$ is the matrix of singular values of the unfolded core tensor $(\mathcal S_{\mathcal{F}})_{(i)}$ along each mode, i.e., $[\sim, \bm \Sigma_{\mathcal{F}_i},\sim] = \texttt{SVD}((\mathcal S_{\mathcal{F}})_{(i)},r'_{\mathcal{F}_i})$, $i=1,2,3$. This quantity measures the relative contribution of the $r$th rank.  
 The rank ($r_{\mathcal{F}_i}$) is adjusted or remains unchanged  to maintain $\epsilon$ within a  desired range of $\epsilon_{l} \leq \epsilon_i \leq \epsilon_{u}$, where $\epsilon_{u}$ and $\epsilon_{l}$ are user-specified upper and lower thresholds.  Rank increase or decrease only requires truncating the SVD of matrix $\bm C_i$ (Eq. \ref{Eq:SVD-Fiber}) at  $r_{\mathcal {F}_i}+1$ or $r_{\mathcal {F}_i}-1$, respectively.  If $r_{\mathcal {F}_i}$ is updated,  then  $r'_{\mathcal {F}_i}=r_{\mathcal{F}_i} + 2$ must be updated accordingly.

\begin{remark}\label{rm:rank}
The size of the intersection tensor along each mode must be greater than or equal to the corresponding Tucker rank, i.e., $r'_{\mathcal {F}_i} \geq r_{\mathcal {F}_i}$. This constraint can be explained by inspecting Eq. \ref{eq:core_lst}. If $r'_{\mathcal {F}_i} = r_{\mathcal {F}_i}$, it is easy to verify that the above least squares problem becomes an interpolation problem, i.e., the  $\hat{\mathcal{F}}(\bm p_1,\bm p_2, \bm p_3) = \mathcal{F}(\bm p_1,\bm p_2, \bm p_3)$. When $r'_{\mathcal {F}_i} > r_{\mathcal {F}_i}$, Eq. \ref{eq:core_lst} amounts to a regression solution, i.e., an overdetermined system of equations. However, if  $r'_{\mathcal {F}_i} < r_{\mathcal {F}_i}$,  Eq. \ref{eq:core_lst} becomes an underdetermined system of equations which can result in poor or unstable solutions. 
\end{remark}

\let\oldnl\nl
\newcommand{\nonl}{\renewcommand{\nl}{\let\nl\oldnl}}

\begin{algorithm}[t]
\begin{footnotesize}
\caption{\small{\texttt{DEIM-FS} Tucker tensor low-rank approximation}}
\label{Algorithm1}
\SetAlgoLined

\KwIn{

$\tilde{\bm U}_{\mathcal{F}}^{(i)}$:  matrix of exact or approximate left singular vectors of $\mathcal{F}_{(i)}$.\\

$r_{\mathcal{F}_i}$:  target Tucker rank.\\

$\mathcal{F}$:  function handle to compute fibers of the target $\mathcal{F}$ }

\KwOut{$\mathcal{S}_{\mathcal{F}},~ \bm U^{(1)}_{\mathcal{F}} ,~ \bm U^{(2)}_{\mathcal{F}},~ \bm U^{(3)}_{\mathcal{F}}$}\vspace{3mm}

$\bm p_i=$ \texttt{DEIM}($\tilde{\bm U}_{\mathcal{F}}^{(i)}$),  \hspace{38mm} $\rhd$ Determine $\bm p_1\in\mathbb{I}^{r'_{\mathcal {F}_1}}$, $\bm p_2\in\mathbb{I}^{r'_{\mathcal {F}_2}}$, $\bm p_3\in\mathbb{I}^{r'_{\mathcal {F}_3}}$ where $r'_{\mathcal {F}_i}={r^{}_{\mathcal{F}}}^{}_i + 2$ \vspace{1mm}

$\bm C_1 = \Bigl(\mathcal{F}(:, \bm p_2, \bm p_3)\Bigr)_{(1)}$, \hspace{25mm} $\rhd$ Calculate $\bm C_1 \in\mathbb{R}^{N_1\times r'_{\mathcal {F}_2}r'_{\mathcal {F}_3}}$, ~$\bm C_2 \in\mathbb{R}^{N_2\times r'_{\mathcal {F}_1}r'_{\mathcal {F}_3}}$, $\bm C_3 \in\mathbb{R}^{N_3\times r'_{\mathcal {F}_1}r'_{\mathcal {F}_2}}$ \vspace{-1mm}

\nonl $\bm C_2 = \Bigl(\mathcal{F}(\bm p_1, :, \bm p_3)\Bigr)_{(2)}$, \vspace{-1mm}

\nonl $\bm C_3 = \Bigl(\mathcal{F}(\bm p_1, \bm p_2, :)\Bigr)_{(3)}$
\vspace{2mm}

$[\bm U^{(i)}_{\mathcal{F}}, \sim , \sim ] = \texttt{SVD}( \bm C_i, r_{\mathcal{F}_i})$ \hspace{17mm}  $\rhd$ Calculate  the left singular vectors of $\bm C_i$ and truncate at rank $r_{\mathcal{F}_i}$ \vspace{2mm}

$\mathcal{W} = \mathcal{F}(\bm p_1, \bm p_2, \bm p_3)$
\hspace{87mm} $\rhd$ Form $\mathcal{W} \in\mathbb{R}^{r'_{\mathcal {F}_1}\times r'_{\mathcal {F}_2}\times r'_{\mathcal {F}_3}}$ \vspace{2mm}

$\mathcal{S}_{\mathcal{F}} = \mathcal{W} \times_1 ~ \bm U^{(1)}_{\mathcal{F}}(\bm p_1, :)^\dagger \times_2 ~ \bm U^{(2)}_{\mathcal{F}}(\bm p_2, :)^\dagger \times_3 ~ \bm U^{(3)}_{\mathcal{F}}(\bm p_3, :)^\dagger$ \hspace{27mm} $\rhd$ Calculate the core tensor ($\mathcal{S}_{\mathcal{F}}$) \vspace{4mm}

\end{footnotesize}
\end{algorithm}

\subsection{Iterative DEIM-FS Tucker tensor approximation}\label{sec:iterative_DEIM-FS}
Step 1 of Algorithm \ref{Algorithm1} requires access to the HOSVD factor matrices or some close approximations to  $\tilde{\bm U}_{\mathcal{F}}^{(i)}$. As we will discuss later, DLRA serves as one such example wherein the factor matrices from the \emph{previous} time step are utilized in the DEIM algorithm to determine which fibers should be sampled at the \emph{current} time step. However, for problems in which these factor matrices are unknown, Algorithm \ref{Algorithm1} can be applied iteratively with minor modifications, as explained below.

In the first iteration, Step 1 is skipped and instead, the fiber indices ($\bm p_i$) are chosen randomly. Then Step 2 is executed. In Step 3, $\tilde{\bm U}_{\mathcal{F}}^{(i)}$ are stored as the first $r'_{\mathcal F_i}$ left singular vectors of matrix $\bm C_i$. Therefore, the first $r_{\mathcal F_i}$ columns of $\tilde{\bm U}_{\mathcal{F}}^{(i)}$  are identical to matrix  $\bm U_{\mathcal{F}}^{(i)}$. However, since $r'_{\mathcal F_i} \geq r_{\mathcal F_i}$ , more SVD columns are stored in matrix $\tilde{\bm U}_{\mathcal{F}}^{(i)}$. Then Steps 4 and 5 are executed. In the second iteration and the iterations after that Steps 1-5 of Algorithm \ref{Algorithm1} are executed with the only modification that $\tilde{\bm U}_{\mathcal{F}}^{(i)}$ are updated in Step 3. The iterations can continue until the singular values of unfolded core tensor (${\mathcal{S}_{\mathcal{F}}}_{(i)}$) converge up to a  threshold value. In particular, the convergence criterion is defined as:
\begin{equation}\label{eq:ErrorProxy}
\dfrac{ \Bigl | \| \boldsymbol \Sigma_{\mathcal{F}_i}^{k} \|_F - \| \boldsymbol \Sigma_{\mathcal{F}_i}^{k-1} \|_F \Bigr | }{\| \boldsymbol \Sigma_{\mathcal{F}_i}^{k} \|_F} < \epsilon, 
\end{equation}
where $\boldsymbol \Sigma_{\mathcal{F}_i}^{K}$ is the matrix of singular values of the $i^{\text{th}}$ unfolding of the core tensor in the $k^{\text{th}}$ iteration and $\epsilon$ is the threshold value. We also note that similar iterative approaches have been used in the past for tensor train cross approximation \cite{ISTtcross}.

\subsection{Computational complexity} \label{sec:ComputationComplx}
In this section, we present the computational complexity of \texttt{DEIM-FS}. In the following analysis,  we use the fact that $r_{\mathcal F_i}$ and $r'_{\mathcal {F}_i}$ are of $\mathcal{O}(r)$. We also assume $N_i \sim \mathcal{O}(N)$. 

 In Step 1, the computational cost of finding the DEIM indices for a $d$-dimensional tensor, scale with $\mathcal{O}(rdN)$. Computing and or storing the fibers in Step 2 scale with $\mathcal{O}(dN r^{d-1})$. 

 The computational cost of performing  SVD is $\mathcal{O}(Nr^{2(d-1)}) +\mathcal{O}(r^{3(d-1)}) $ or $\mathcal{O}(N^2r^{d-1}) +\mathcal{O}(N^3) $, whichever is smaller. One can also use randomized SVD to obtain  $\bm U^{(i)}_{\mathcal{F}}$ in order to reduce the computational cost of performing SVD.  Randomized SVD algorithms can be particularly effective since the $\bm C_i$ matrices have a large number of columns ($\mathcal{O}(r^{d-1})$) and only a small ($r$) left singular vectors need to be computed accurately. We have not used randomized SVD algorithms in any of the test cases in this paper.

Forming  $\mathcal W$ does not require extra calculation as the values of the intersection tensor are already calculated in Step 2 and can be extracted from any of the $\bm C_i$ matrices. The computational cost of computing the core tensor is $\mathcal{O}(dr^{d+1})$. 

In summary, Steps 2 and 3 are the two most computationally expensive parts of the \texttt{DEIM-FS} algorithm. As mentioned above, in Step 3, randomized SVD algorithms can be utilized to further reduce the cost. On the other hand, Step 2 requires accessing elements of tensor $\mathcal{F}$. For problems where accessing any element of tensor $\mathcal{F}$ requires additional computation, the computational cost of computing the fibers can dominate the overall cost. An example of this type of problem is DLRA, where computing any element of the right-hand side tensor $\mathcal{F}=\mathscr{F}(\hat{\mathcal{V}})$ requires applying a nonlinear discrete differential operator on tensor $\hat{\mathcal{V}}$. 

\subsection{DEIM fiber sampling for dynamical low-rank approximation}  \label{sec:DEIM-FS DLRA}
The \texttt{DEIM-FS} algorithm can be employed to construct a low-rank Tucker tensor approximation of $\mathscr{F}(\hat{\mathcal{V}})$ involved in the DLRA evolution equations. As demonstrated, this approach leads to significant computational savings when solving Eqs. (\ref{eq:Sevolution}) and (\ref{eq:Uevolution}).
Step 1 of Algorithm \ref{Algorithm1} requires $\tilde{\bm U}^{(i)}_{\mathcal{F}}$, representing the exact or approximate left singular vectors of matrices $\mathcal{F}_{(i)}$. To obtain these vectors, we utilize $\bm U^{(i)}_{\mathcal{F}}$ from the previous time step, which has already been calculated using the \texttt{DEIM-FS} algorithm. It is important to note that the $\tilde{\bm U}^{(i)}_{\mathcal{F}}$ matrices are necessary solely for the DEIM algorithm to compute the fiber indices, while the actual computation of the fibers is carried out for $\mathscr{F}(\hat{\mathcal{V}})$ at the current time step.

In practice, utilizing $\tilde{\bm U}_{\mathcal{F}}^{(i)}$ from the previous time step yields excellent performance. The difference in accuracy compared to cases where Algorithm \ref{Algorithm1} is used iteratively, as detailed in Section \ref{sec:iterative_DEIM-FS}, is negligible. The same approach was used in previous studies for low-rank approximation of matrix differential equations \cite{MNAdaptive,DPNFB23}. At $t=0$, if no good approximation for $\tilde{\bm U}^{(i)}_{\mathcal{F}}$ exists, \texttt{DEIM-FS (iterative)} may be used. 

In the DLRA equations, the fibers of tensor $\mathcal{F}$ are calculated by evaluating the function $\mathcal{F} = \mathscr{F}(\hat{\mathcal{V}})$. For the TDEs obtained from discretizing a PDE, the function $\mathscr{F}( \sim )$ involves discrete differential operators. 
 Since the derivative calculation requires adjacent points of the DEIM-selected points, the adjacent points must be determined in Step 2 of  Algorithm \ref{Algorithm1}. This step depends on the numerical scheme used for the spatial
discretization. For instance, when employing the spectral element method, determining the derivative at a specific spatial point requires access to the values of other points within the same element \cite{MNAdaptive}. As a result, calculating any fiber of $\mathcal{F}$ requires access to the values to additional (adjacent) fibers of $\hat{\mathcal{V}}$. We denote the additional fiber indices with $\bm p_{a_i}$.  For example, for a three-dimensional tensor, $\bm C_1$ is calculated as:
\begin{equation}\label{eq:C1}
 \mathcal{F}(:, \bm p_2, \bm p_3) = \mathscr{F} \Bigl( \mathcal{S} \times_1 ~ \bm U^{(1)} \times_2 ~ \bm U^{(2)}([\bm p_2, \bm p_{a_2}],:) \times_3 ~ \bm U^{(3)}([\bm p_3, \bm p_{a_3}]  ,:) \Bigr).
\end{equation}

Applying this methodology to DLRA enables a significant reduction in the computational cost and memory. In these equations, $\mathcal{F}$ is $N^d$ tensor, while the memory advantage of \texttt{DEIM-FS} anbles storing $\mathcal{F}$ in the Tucker compressed form. Accordingly, Eqs. (\ref{eq:Sevolution}) and (\ref{eq:Uevolution}) can be rewritten as:
\begin{subequations}
\begin{gather}
  \Dot{\mathcal{S}}  = \mathcal{S}_{\mathcal{F}} ~ \bigtimes_{i=1}^d ~ \bigl( {\bm U^{(i)}}^T \bm U^{(i)}_{\mathcal{F}} \bigr) , 
  \label{eq:SevolutionReduced} \\
  \Dot{\bm U}^{(i)}  = \Bigl(\bm I - \bm U^{(i)} \bm U^{(i)^T} \Bigr) \bm U^{(i)}_{\mathcal{F}} ~ \Bigl[ \mathcal{S}_{\mathcal{F}} \bigtimes_{k \neq i} \bigl( {\bm U^{(k)}}^T \bm U^{(k)}_{\mathcal{F}} \bigr) \Bigr]_{(i)} ~ \mathcal{S}_{(i)}^\dag , 
  \label{eq:UevolutionReduced}
\end{gather}
\end{subequations}\\
where $\mathcal{S}_{\mathcal{F}}$ and $\bm U^{(i)}_{\mathcal{F}}$ are the outputs of the Algorithm \ref{Algorithm1}. Another advantage of employing \texttt{DEIM-FS} for DLRA is its simplification of the implementation of DLRA equations. This is due to the tensor $\mathcal{F}$ being approximated in a black-box fashion and Eqs. (\ref{eq:SevolutionReduced}) and (\ref{eq:UevolutionReduced}) are agnostic to the type of TDE being solved by DLRA.  Hereinafter, Eqs. (\ref{eq:SevolutionReduced}) and (\ref{eq:UevolutionReduced}) are referred to as \texttt{DLRA-DEIM-FS}.

\subsection{Comparison to an existing fiber sampling algorithm}
\label{sec:ComparetoFSTD}
We compare our proposed algorithm with the FSTD algorithm \cite{CCGeneralizing}. While FSTD has some similarities to \texttt{DEIM-FS}   it also has some key differences with the presented algorithm. We briefly review the FSTD algorithm here. In FSTD, the cross Tucker model of  $\mathcal{F}$ is given by:
\begin{equation}\label{eq:FSTD}
\mathcal{F} \approx \mathcal{W} \times_1 \bm C_1  \mathcal{W}_{(1)}^{\dagger} \times_2 \bm C_2  \mathcal{W}^\dag_{(2)} \times_3 \bm C_3  \mathcal{W}^\dag_{(3)},
\end{equation}
where $\bm C_1 \in \mathbb{R}^{N_1 \times r'_{\mathcal {F}_2}r'_{\mathcal {F}_3}}$, $\bm C_2 \in \mathbb{R}^{N_2 \times r'_{\mathcal {F}_1}r'_{\mathcal {F}_3}}$, $\bm C_3 \in \mathbb{R}^{N_3 \times r'_{\mathcal {F}_1}r'_{\mathcal {F}_2}}$ are the unfolded $\mathcal{C}_1, \mathcal{C}_2, \mathcal{C}_3$ fiber sub-tensors shown in Figure \ref{fig:DEIM-FS}, $\mathcal{W} \in \mathbb{R}^{r'_{\mathcal {F}_1} \times r'_{\mathcal {F}_2} \times r'_{\mathcal {F}_3}}$ is the intersection tensor. According to the FSTD algorithm presented in \cite{CCGeneralizing}, at the first step, the index of the first fiber is initialized. Then, the indices of the other fibers are determined based on a deterministic greedy algorithm.

We compare FSTD to \texttt{DEIM-FS} in terms of the information these algorithms require about the tensors, their accuracy, and the computational cost. The FSTD algorithm does not require any information about $\mathcal{F}$, whereas  \texttt{DEIM-FS} requires the exact or approximate left singular vectors of all unfoldings of $\mathcal{F}$.  The \texttt{DEIM-FS (iterative)} algorithm, on the other hand, does not require left singular vectors of $\mathcal{F}$ unfolding. However, computing \texttt{DEIM-FS (iterative)} requires iterations and it is more expensive to compute in comparison to FSTD. 

The   \texttt{DEIM-FS}  Tucker tensor models appear to be more accurate than FSTD Tucker tensor models of the same rank. This is due to several factors: (i) The FSTD becomes ill-conditioned as the rank increases. The same issue also exists for matrix CUR decompositions, where the intersection matrix becomes singular as rank increases \cite{ISTtcross}. This issue is mitigated in matrix CUR by performing QR decomposition of selected columns or rows.  The issue of ill-conditioning is also mitigated in \texttt{DEIM-FS}, since the factor matrices are obtained by performing SVD of the selected fibers. The SVD, similar to QR, results in a set of orthonormal modes, and the matrices $\bm U^{(i)}_{\mathcal F}(\bm p_i,:)$ are well-conditioned matrices. In fact, the DEIM algorithm is designed to maintain  $\| \bm U^{(i)}_{\mathcal F}(\bm p_i,:)^{-1} \|_2$ as small as possible. (ii) The choice of initial fibers in FSTD can have a significant impact on the accuracy of the resulting Tucker model. We show this effect in our numerical examples. We also show that \texttt{DEIM-FS (iterative)} is much less sensitive to the initial choice of fibers.  

From the computational cost point of view, both FSTD and \texttt{DEIM-FS} require access to $\mathcal{O}(dr^{d-1}N)$ number of elements of $\mathcal{F}$. They also have the same memory requirements.    However, \texttt{DEIM-FS} requires more flops due to the computation of SVD of $\bm C_i$ matrices. The FSTD algorithm offers the advantage of storing the actual fibers of the target tensor, thereby inheriting the structure of the target tensor. For instance, in the case of sparse tensors, storing FSTD factor matrices in a sparse form can potentially achieve a very high compression ratio.

\section{Demonstration}
\label{sec:Demo}

\subsection{Toy Examples} 
\label{sec:ToyEx}
As our first example, we consider two three-dimensional functions as shown below:
\begin{subequations}
\begin{align}
  \mathscr{F}_1(x_1,x_2,x_3)  = e^{-(x_1 ~ x_2 ~ x_3)^2} \hspace{3.4cm}  x_1, x_2, x_3 \in [-1,1] ,
  \label{eq:toy1} \\
  \mathscr{F}_2(x_1,x_2,x_3)  = \frac{1} {\bigl(x_1^b + x_2^b + x_3^b\bigr)^{1/b}} \hspace{1cm} x_1, x_3 \in [1,300]  ~~~~ x_2 \in[1,400] .
  \label{eq:toy2}
\end{align}
\end{subequations}

The tensor $\mathcal{F}_1$ is obtained by evaluation $\mathscr{F}_1(x_1,x_2,x_3)$ at 100 equally spaced elements of $x_1$, $x_2$, and $x_3$ in their respective domains. Therefore, $\mathcal{F}_1 \in \mathbb{R}^{100 \times 100 \times 100}$. Similarly,  $\mathcal{F}_2 \in \mathbb{R}^{300 \times 400 \times 300}$ is obtained by evaluating $\mathscr{F}_2(x_1,x_2,x_3)$ on a uniform grid in each direction. Two choices of $b=3$ and $b=5$ are considered for $\mathcal{F}_2$. In all demonstrations of \texttt{DEIM-FS}, the exact left singular vectors 
 of the unfolding of the target tensors are used in the DEIM algorithm.   To compare the accuracy and efficiency of \texttt{DEIM-FS} against HOSVD and FSTD, \cite{CCGeneralizing}, $\mathcal{F}_1$ and $\mathcal{F}_2$ are approximated using the three mentioned algorithms. Then the error between the approximated tensors and the actual tensors is calculated. Denoting $\hat{\mathcal{F}}_1$ and  $\hat{\mathcal{F}}_2$ as the low-rank Tucker approximation of $\mathcal{F}_1$ and $\mathcal{F}_2$, respectively, the error is  defined as:
\begin{equation}\label{eq:Error1}
 \mathcal{E}  = \| \hat{\mathcal{F}} - \mathcal{F} \|_F .
\end{equation}
For \texttt{DEIM-FS}, we consider ${r_{\mathcal{F}}}^{}_1 = {r_{\mathcal{F}}}^{}_2 = {r_{\mathcal{F}}}^{}_3 = {r_{\mathcal{F}}}$. Hence, $r'_{\mathcal {F}_1} = r'_{\mathcal {F}_2} = r'_{\mathcal {F}_3} = r'_{\mathcal {F}}$. As mentioned in Remark 1, $r'_{\mathcal {F}} \geq r_{\mathcal {F}}$.   In Figure 2a, we increase $r'_{\mathcal {F}}$ for a fixed ${r_{\mathcal{F}}}$ to study the effect of increasing $r'_{\mathcal {F}}$. The results of Figure 2a indicate that increasing $r'_{\mathcal {F}}$ beyond $r_{\mathcal{F}} + 2$ results in a negligible reduction in error. We use $r'_{\mathcal{F}} = r_{\mathcal{F}} + 2$ for rest of examples in this paper. 

Figures 2b and 2c show the low-rank error versus rank for $\mathcal{F}_1$ and $\mathcal{F}_2$, respectively. As it can be seen, for both tensors, the error of \texttt{DEIM-FS} method closely follows the HOSVD error. However, the FSTD error is always greater than the \texttt{DEIM-FS} error. Moreover,  the FSTD error either does not decrease (see Figure 2c) or even increases (see Figure 2b)  as rank increases. This is due to the issue of ill-conditioning as discussed in Section \ref{sec:ComparetoFSTD}.

\begin{figure}[htbp]
     \centering
     \includegraphics[scale=1]{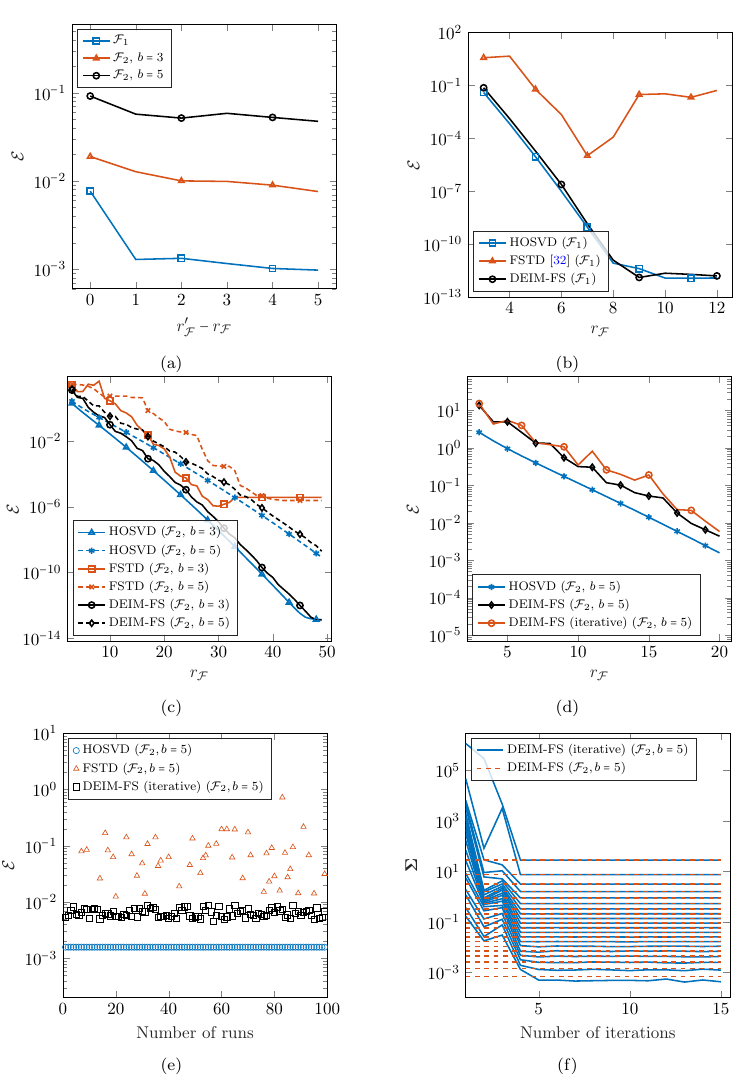}
     \caption{Toy Example: Comparison of the FSTD, \texttt{DEIM-FS},  \texttt{DEIM-FS (iterative)}, and HOSVD: (a) approximation errors of $\mathcal{F}_1$ and $\mathcal{F}_2$ versus $r'_{\mathcal{F}}$ for a fixed rank of $r_{\mathcal{F}}$; (b) approximation error of $\mathcal{F}_1$ versus rank; (c) approximation error of $\mathcal{F}_2, (b=3,5)$ versus rank; (d) approximation error of $\mathcal{F}_2, (b=5)$ versus rank; (e) effect of fiber initialization on the approximation error of $\mathcal{F}_2, (b=5)$ for a fixed rank of $r_{\mathcal{F}}=20$ for 100 different initializations; (f) convergence of the singular values of the core tensor ($\mathcal{S}_{\mathcal{F}}$) versus iterations for \texttt{DEIM-FS (iterative)} algorithm.}
     \label{fig:CompareMethods}
\end{figure}

As discussed in Section \ref{sec:iterative_DEIM-FS}, in cases where $\tilde{\bm U}_{\mathcal{F}}^{(i)}$ is not available, the \texttt{DEIM-FS (iterative)} approach can be employed. Figure 2d illustrates that the error obtained by the \texttt{DEIM-FS (iterative)} algorithm closely follows the errors of the HOSVD and \texttt{DEIM-FS} methods. In the \texttt{DEIM-FS (iterative)} algorithm, the initial fibers are randomly selected, while in \texttt{DEIM-FS}, the exact left singular vectors are utilized.

In both the FSTD and \texttt{DEIM-FS (iterative)} algorithms, the initial fibers are selected randomly. Figure 2e examines the effect of different initializations on accuracy. The errors of Tucker models, for $\mathcal{F}_2$ resulting from 100 random initializations of both FSTD and \texttt{DEIM-FS (iterative)}, are depicted in Figure 2e. It is noticeable that the variance of errors obtained by FSTD is significantly larger than those obtained by \texttt{DEIM-FS (iterative)}. Interestingly, for some fiber initializations in FSTD, the fiber indices remain fixed on the initial fibers and do not update, leading to very large errors ($\mathcal E \approx 30$). These cases are not displayed in Figure 2e because we believe a minor fix could resolve this issue. In summary, our observation reveals that the error in FSTD is highly sensitive to fiber initializations, while \texttt{DEIM-FS (iterative)} displays a significantly smaller variance in error.

For the 100 initializations presented in Figure 2e, the average number of iterations required for convergence by \texttt{DEIM-FS (iterative)} is 5.48, with a standard deviation of 1.24. Figure 2f illustrates the convergence of the singular values of the core tensor ($\mathcal{S}_{\mathcal{F}}$) determined by the  \texttt{DEIM-FS (iterative)} algorithm.  

\subsection{DLRA for four-dimensional Fokker Planck  equation}
\label{sec:FPexample}
For the second test case, a four-dimensional Fokker Planck (FP) equation is considered. The physics under consideration is governed by the Langevin form of stochastic differential equations (SDE).
\begin{equation}\label{eq:Langevin}
 d \bm x = \bm f(\bm x,t)dt + \bm S (\bm x,t) d\bm w ,
\end{equation}
where $\bm x=[x_1, x_2, x_3, x_4]$ is the stochastic transport variable, $\bm f(\bm x,t)$ is the drift function, $\bm S(\bm x,t)$ is the diffusion matrix and $\bm w$ denotes the Wiener-Levy process \cite{SKSecond}. In this equation, $\bm x$ and $\bm f(\bm x,t)$ are $d$-dimensional vectors, $\bm S(\bm x,t)$ is a $d \times m$ matrix, and $\bm w$ is an $m$-dimensional vector. The transitional probability density function (PDF) of the stochastic variable is governed by its FP equation:
\begin{equation}\label{eq:FokkerPlanck}
 \frac{\partial p(\bm x,t)}{\partial t} = - \mathlarger{\sum}_{i=1}^{d} \frac{\partial}{\partial x_i} \bigl[  f_i(\bm x,t) p(\bm x,t) \bigr] + \frac{1}{2} \mathlarger{\sum}_{i=1}^{d} \mathlarger{\sum}_{j=1}^{d} \frac{\partial^2}{\partial x_i \partial x_j} \bigl[ D_{ij}(x,t) p(x,t) \bigr] ,
\end{equation}
where $\bm D=[D_{ij}]$ is a $d \times d$ matrix and $\bm D = \bm S \bm S ^T$. We consider a 4-dimensional super-symmetric PDF and we choose $r = r_1 = r_2 = r_3 = r_4$ and ${r_{\mathcal{F}}}_1 = {r_{\mathcal{F}}}_2 = {r_{\mathcal{F}}}_3 = {r_{\mathcal{F}}}_4 = {r_{\mathcal{F}}}$. We consider homogeneous Dirichlet boundary condition and a normal distribution as the initial condition as: $p_0(\bm x) = \mathcal{N} \bigl(\boldsymbol \mu_0,~ \bm C_0 \bigr)$, 
where $\mathcal{N}$ is the normal distribution and $\boldsymbol \mu_0$ is the mean set to $\boldsymbol \mu_0 
            = \begin{bmatrix}  1.5 & 0.6 & -0.3 & -1.2
            \end{bmatrix}^T$, and  $\bm C_0$ is the covariance matrix, where $\bm C_{0_{ii}} = 1$, $i=1,\dots, 4$ and $\bm C_{0_{ij}} = 0.5$,  $i,j=1,\dots, 4$ and $i\neq j$.
 We also consider $f_i = -\alpha x_i~,~i=1,2,3,4$ and $\alpha = 0.75$. The solution $p(x,t)$ remains Gaussian when $D_{ij}$ is a constant. The analytical expressions are obtained for the moments over time as:
\begin{equation}\label{eq:MeanOverTime}
 \boldsymbol \mu(t) = \boldsymbol \mu_0 \exp(-\alpha t) , 
\end{equation}
\begin{equation}\label{eq:CovOverTime}
 \bm C(t) = \frac{\bm D}{2\alpha} + \bigl( \bm C_0 - \frac{\bm D}{2\alpha} \bigr)  \exp(-2\alpha t).
\end{equation}
For spatial discretization, the third-order spectral method is used.  Each domain is discretized via $N = N_1 = N_2 = N_3 = N_4 = 61$ points, which includes 20 elements with a second-order order polynomial approximation within each of the elements. For time advancement, the fourth-order Runge-Kutta (RK4) method is used. The domain under consideration is $x \in [-6,6]^4$ within the time interval $t \in [0,8]$. The time advancement of  $\Delta t = 2 \times 10^{-3}$ is chosen. The resulting TDE is solved using \texttt{DLRA-DEIM-FS} (Eqs. \ref{eq:SevolutionReduced} and \ref{eq:UevolutionReduced}) with $r=5$. Although the algorithm can be adaptive, we solve this example with fixed ${r_{\mathcal{F}}}=5$ and $r'_{\mathcal{F}}=7$. Figure 3a shows the temporal evolution of the singular values obtained from \texttt{DLRA-DEIM-FS} and analytical solution. For the rest of the paper, the error is defined as the relative error as:
\begin{equation}\label{eq:Error2}
 \mathcal{E}(t)  = \dfrac{\| \hat{\mathcal{V}}(t) - \mathcal{V}(t) \|_F}{\| \mathcal{V}(t) \|_F}.
\end{equation}

\begin{figure}[t]
     \centering
     \includegraphics[scale=1]{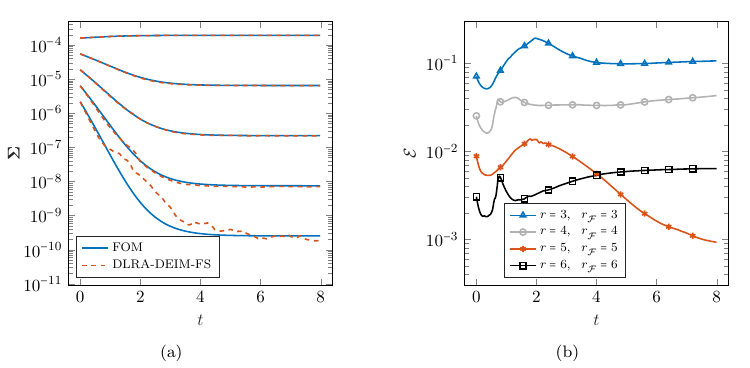}
     \caption{Fokker Planck Equation: (a) first 5 singular values of analytical solution and \texttt{DLRA-DEIM-FS}; (b) relative error evolution for different $r$ and $r^{}_{\mathcal{F}}$.}
     \label{fig:Ex2Results}
\end{figure}

We conduct a convergence study by solving the DLRA equations with various $r$ and $r_{\mathcal{F}}$ values. The temporal evolution of the error is depicted in Figure 3b. It is observed that the model with $r=5$ and $r_{\mathcal{F}}=5$ yields the best error. The statistical moments obtained by solving the FP equation with \texttt{DLRA-DEIM-FS} can be compared to those obtained analytically (refer to Eqs. (\ref{eq:MeanOverTime}) and (\ref{eq:CovOverTime})). Table \ref{table:Ex2Moments} presents the comparison of the mean of the variables and the covariance matrices at the final time step $(t=8)$.

\begin{table}[t]
\vspace{-4mm}
\caption{Comparison of the moments of the PDF at $t=8$ with $r=5$ and $r^{}_{\mathcal{F}} = 5$}
\vspace{-3mm}
\fontsize{8.5pt}{8pt}\selectfont
\centering
        \begin{tabular}{  c  c  c  }
        \hline
             & &  \\
             & Mean & Covariance \\
            \hline
            & &  \\
            Analytical ~~~ & $\begin{bmatrix}
            0.0000 & 0.0000 & 0.0000 & 0.0000
            \end{bmatrix}$ & $\begin{bmatrix} 
            0.6675 & 0.0272 & 0.0272 & 0.0272 \\
            0.0272 & 0.6675 & 0.0272 & 0.0272 \\
            0.0272 & 0.0272 & 0.6675 & 0.0272 \\
            0.0272 & 0.0272 & 0.0272 & 0.6675 \\
            \end{bmatrix}$ \\ 

             & &  \\
            \hline
             & &  \\
             
            DLRA-DEIM-FS ~~~ &  $\begin{bmatrix}
            0.0017 & 0.0041 & -0.0007 & -0.0028
            \end{bmatrix}$ & $\begin{bmatrix} 
            0.6674 & 0.0272 & 0.0270 & 0.0271 \\
            0.0272 & 0.6673 & 0.0271 & 0.0272 \\
            0.0270 & 0.0271 & 0.6675 & 0.0270 \\
            0.0271 & 0.0272 & 0.0270 & 0.6674 \\
            \end{bmatrix}$ \\
             & &  \\
        \hline
        \end{tabular}
\label{table:Ex2Moments}
\end{table}

Although this example is not a nonlinear TDE, there are still advantages to using \texttt{DLRA-DEIM-FS} to solve this TDE rather than employing standard DLRA techniques. To illustrate this, consider the right-hand side of the FP equation, which involves the summation of 20 terms. Implementing standard DLRA requires a highly intrusive and meticulous treatment of these terms to avoid storing the full-dimensional right-hand side tensor. However, the \texttt{DLRA-DEIM-FS} algorithm remains agnostic to the nature of the TDE, involving the evaluation of sparse fibers of the right-hand side in a black-box fashion. Therefore, \texttt{DLRA-DEIM-FS} is significantly easier to implement compared to the standard practice for solving DLRA equations, even for linear TDEs.

Solving the FOM for this example requires storing $N^4 = 13,845,841$ floating-point numbers for the right-hand side of the TDE. Using \texttt{DLRA-DEIM-FS} with $r_{\mathcal{F}}=5$ requires storing $d \times N \times {r'_{\mathcal{F}}}^{d-1} + {r'_{\mathcal{F}}}^4 = 86,093$ floating-point numbers, resulting in a memory compression ratio of $13,845,841 / 86,093 = 160.8$.

\subsection{DLRA for four-dimensional nonlinear advection equation}
For the final demonstration, we consider a four-dimensional nonlinear advection given below:
\begin{equation}\label{eq:NonlinearPDE}
 \frac{\partial v(\bm x,t)}{\partial t} = -\bm b \cdot \nabla v(\bm x,t) + s(v(\bm x,t)), \quad \bm x \in [-5,5]^4,
\end{equation} 
where $\bm b = [-\sin(t),~ \cos(t),~ -\sin(\pi + t),~ \cos(\pi + t)]$ is the advection velocity and $s(v)=- 0.1v/(1+v^2)$ is the nonlinear source term.  The boundary condition is periodic and the initial condition is considered as  $v_0(\bm x) = f(\bm x) + g(\bm x)$, where 
\begin{equation*}
 f(\bm x) = e^{-(x_1-\frac{1}{2})^2}  e^{-(x_1+\frac{x_2}{2}-\frac{1}{2})^2} e^{-(x_3-\frac{1}{2})^2} e^{-(x_3+\frac{x_4}{2}-\frac{1}{2})^2},
\end{equation*}
\begin{equation*}
 g(\bm x) =  e^{-(x_1+\frac{1}{2})^2} e^{-(x_2+\frac{1}{2})^2}  e^{-(x_3+\frac{1}{2})^2}  e^{-(x_4+\frac{1}{2})^2}.
\end{equation*}
 The second-order spectral method is used for spatial discretization ($\bm x$), and the fourth-order Runge-Kutta (RK4) method is employed for temporal integration within the time interval $t \in [0,4]$. The time step is set to $\Delta t = 2 \times 10^{-3}$. Each domain is discretized with $N = N_1 = N_2 = N_3 = N_4 = 85$ points, comprising  42 elements with a second-order polynomial approximation within each element. The rank-adaptive \texttt{DLRA-DEIM-FS} is utilized for approximating the right-hand side of the resulting TDE. It is important to note that the rank ($r$) of the solution $\hat{\mathcal{V}}(t)$ is not adaptive. The relative error with respect to the norm (Eq. (\ref{eq:Error2})) is used for error calculation. As the analytical solution is unavailable in this case, the FOM is solved using $N=85$ and $\Delta t = 2 \times 10^{-3}$. The solution obtained from FOM is considered the ground truth.

 Three distinct errors contribute to the overall error while using \texttt{DLRA-DEIM-FS}: the \texttt{DEIM-FS} error for low-rank approximation of $\mathscr{F}(\hat{\mathcal{V}})$, controlled by varying threshold values; the DLRA error for approximating $\hat{\mathcal{V}}(t)$, controlled by $r$; and the temporal integration error, influenced by adjusting $\Delta t$.

First, the TDE is solved with a fixed $r=6$ and various rank-adaptivity thresholds: $\epsilon_u={10^{-2}, 10^{-4}, 10^{-6}}$ and $\epsilon_l={10^{-3}, 10^{-5}, 10^{-7}}$. Figure 4a indicates that for $(\epsilon_l=10^{-6},\epsilon_u=10^{-7})$, the DLRA error does not decrease. This suggests that the error in approximating the right-hand side tensor is not dominant, and the error reaches a saturation point due to the DLRA low-rank approximation error in $\hat{\mathcal{V}}(t)$. This also demonstrates that in DLRA it is not necessary to compute  $\mathscr{F}(\hat{\mathcal{V}})$ exactly.  

Figure 4b illustrates the evolution of $r_{\mathcal{F}}$ values associated with different $\epsilon_l$ and $\epsilon_u$ threshold values. As expected, lower thresholds correspond to higher ranks. In Figure 4c, we examine the effect of the rank $\hat{\mathcal{V}}(t)$ on the total error. In this case, the threshold values for approximating $\mathscr{F}(\hat{\mathcal{V}})$ are fixed at $(\epsilon_l=10^{-5},\epsilon_u=10^{-4})$. The results demonstrate that increasing $r$ enhances accuracy. Another study is also conducted to consider the effect of the time step $(\Delta t)$ on the error of the \texttt{DLRA-DEIM-FS} method. To this end, the DLRA rank  and the thresholds are fixed: $r=11$ and  $(\epsilon_l = 10^{-5},\epsilon_u = 10^{-4})$. The TDE is then solved for three different values of $\Delta t=\{2\times 10^{-3}, 4\times 10^{-3},6\times 10^{-3} \}$. Note that the $\Delta t$ of the FOM (the ground truth) was not changed and it was fixed to be $2 \times 10^{-3}$. As shown in Figure 4d the effect of different time steps is considerable until $t \approx 2$. Afterward, the low-rank approximation errors of $\hat{\mathcal{V}}(t)$ and $\mathscr{F}(\hat{\mathcal{V}})$ dominates, and decreasing $\Delta t$ does not change the error considerably.

Figure 5a presents a comparison between the first 11 singular values obtained using \texttt{DLRA-DEIM-FS} and those of FOM, while Figure 5b shows the corresponding evolution of $r_{\mathcal{F}}$ for this problem.

\begin{figure}[t]
     \centering
     \includegraphics[scale=1]{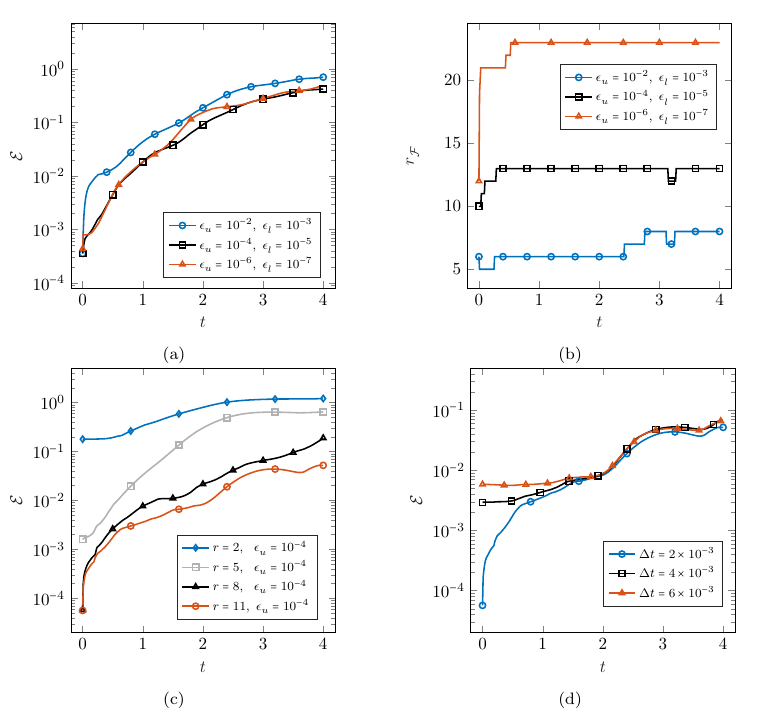}
     \caption{Four-dimensional nonlinear advection equation (Example 3): (a) Relative error evolution with fixed $r$ for different $\epsilon_l$ and $\epsilon_u$; (b) Evolution of $r^{}_{\mathcal{F}}$ associated with the plots of Figure 4a; (c) Relative error evolution for different $r$ with fixed $\epsilon_l = 10^{-5}$ and $\epsilon_u = 10^{-4}$; (d) Relative error evolution for different $\Delta t$ with fixed $r$, $\epsilon_l$, and $\epsilon_u$.}
     \label{fig:Ex3Results1}
\end{figure}

\begin{figure}[t]
     \centering
     \includegraphics[scale=1]{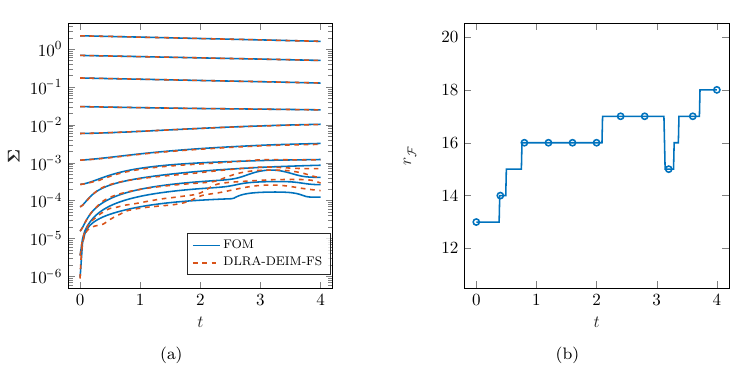}
     \caption{Four-dimensional nonlinear advection equation (Example 3): (a) First 11 singular values of FOM and \texttt{TDB-DEIM-FS}; (b) Evolution of $r^{}_{\mathcal{F}}$ associated with the solved system in Figure 5a.}
     \label{fig:Ex3Results2}
\end{figure}

We also compute the following marginal function:
\begin{equation}\label{eq:MarginPlots}
 \overline{v}(x_3, x_4, t) = \int_{-5}^{5} \int_{-5}^{5} {v(\bm x,t)}^2 ~dx_1 dx_2 ,
\end{equation}
using both \texttt{DLRA-DEIM-FS} and  FOM. 

\begin{figure}[t]
     \centering
     \includegraphics[scale=0.35]{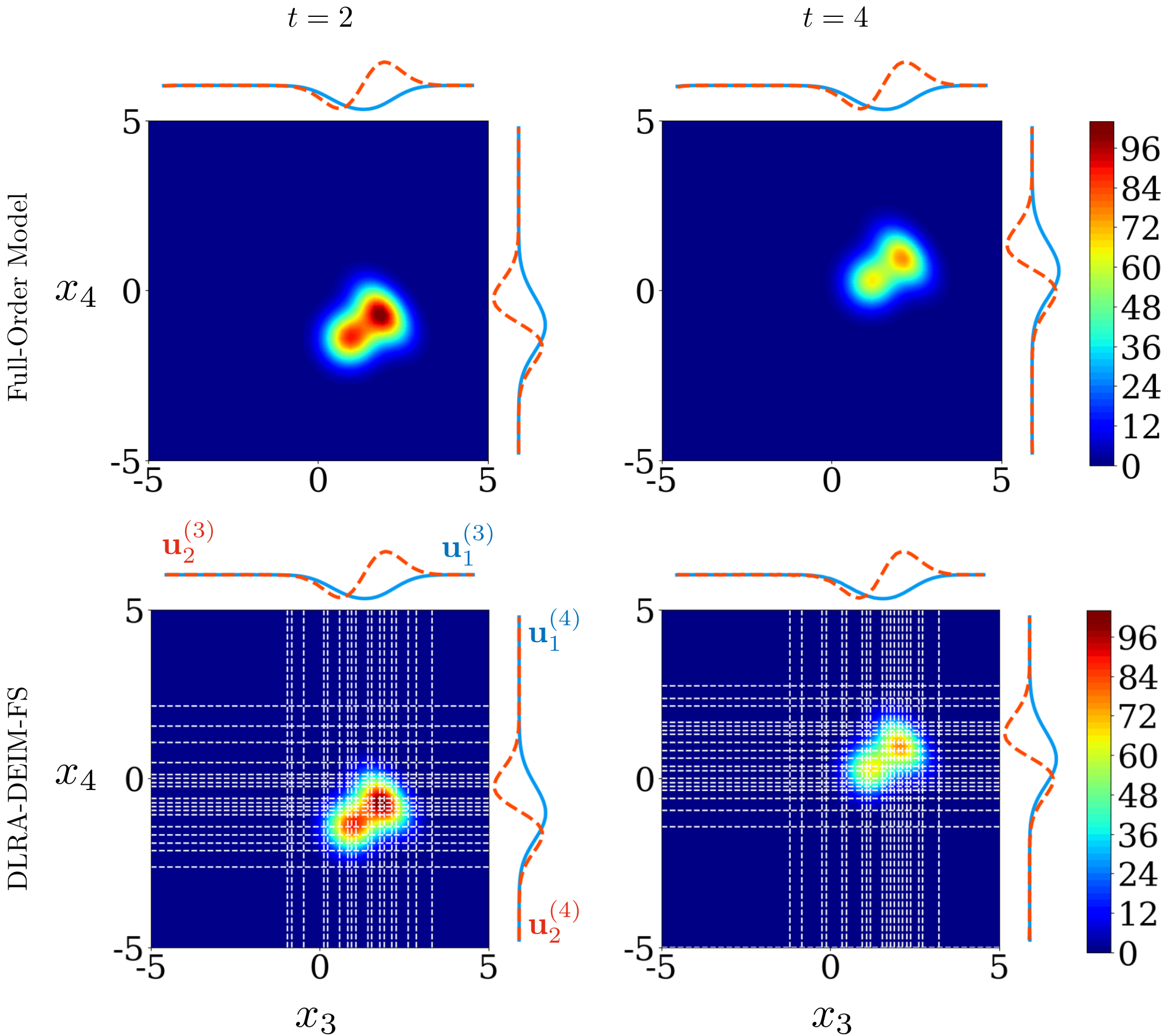}
     \caption{Four-dimensional advection-reaction equation (Example 3): Top row: FOM; Bottom row: \texttt{DLRA-DEIM-FS}.  The DEIM sampled fibers along each mode are shown with white dashed lines. On the top of each panel the first two $x_3$-bases (the first two columns of $\bm U^{(3)}$ ) are shown. The first basis is shown in solid blue line and the second basis is shown in dashed red line. On the right side of each panel. }
     \label{fig:Ex3Marginals}
\end{figure}

Figure \ref{fig:Ex3Marginals} shows contour plots of $\overline{v}(x_3,x_4)$ obtained from FOM (top row) and \texttt{DLRA-DEIM-FS} (bottom row) at two different time instances.   The first two columns of the $x_3$ and $x_4$ factor matrices are shown: $\{ \bm u_1^{(3)}, \bm u_2^{(3)}  \}$ are shown on the top side of each panel and $\{ \bm u_1^{(4)}, \bm u_2^{(4)}  \}$ are shown on the right side of each panel. The first basis is displayed using solid blue color, while the second basis is shown with red dashed lines.  The factor matrices of the FOM are obtained by performing HOSVD on the FOM solution tensor.  The selected fibers along the third and fourth modes of the tensor are also shown with dashed lines.   A good agreement between the \texttt{DLRA-DEIM-FS} and FOM contour plots is observed. Also, it can be seen that the first two $\bm x_3$ and $\bm x_4$ bases are localized and advect with the solution. Moreover, the selected fibers are concentrated around the localized solution. 

The PDE defined by Eq. \ref{eq:NonlinearPDE} contains non-polynomial nonlinearity. Consequently, in standard DLRA implementation, the explicit formation of the right-hand-side tensor is inevitable because $\mathcal F =\mathscr{F}(\hat{\mathcal V})$ is full rank despite $\hat{\mathcal V}$ having a low rank. The explicit formation of $\mathcal F$ incurs substantial memory and computational costs, scaling at least as $\mathcal{O}(N^4)$, both in terms of memory allocation and floating-point operations. The \texttt{DLRA-DEIM-FS} algorithm avoids this cost because the full-rank tensor is never formed.  Instead, a cross Tucker tensor approximation of $\mathcal F$ is constructed on the fly. The computational cost of computing  $\mathcal F$ using \texttt{DLRA-DEIM-FS} scales linearly with $N$.  In Figure 7a the computational cost of solving DLRA evolution equations versus $N$ using standard DLRA and \texttt{DLRA-DEIM-FS} are shown.

 Figure 7b displays the required memory relative to FOM for storing $\mathcal{F} = \mathscr{F}(\hat{\mathcal V})$ in Examples 2 and 3. The results demonstrate a significant reduction in memory for both of these examples.

\begin{figure}[t]
     \centering
     \includegraphics[scale=1]{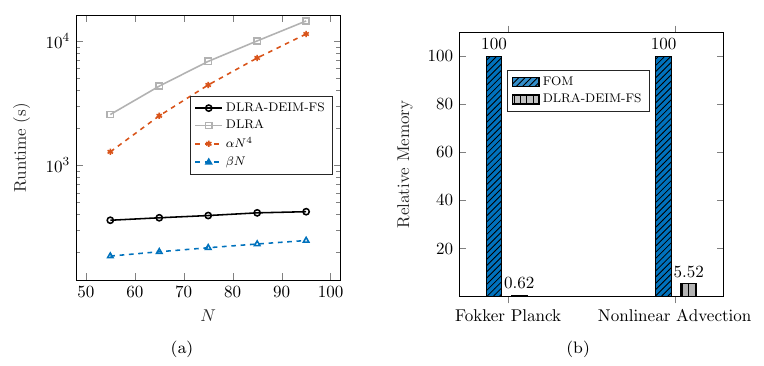}
     \caption{(a) Computational time vs $N$ (number of discretization points) for 4D nonlinear advection equation Four-dimensional advection-reaction equation (Example 3); (b) Comparison of the relative memory for storing the RHS for 4D Fokker Planck and 4D nonlinear advection equations. }
     \label{fig:SpeedUps}
\end{figure}

\section{Conclusions}
\label{sec:conclusions}
Tucker tensor low-rank approximation is a building block of many important tensor low-rank approximation algorithms. Determining the optimal (lowest error) Tucker tensor model lacks a known closed solution. Instead,  HOSVD is commonly used, which yields near-optimal Tucker tensor models. However, HOSVD requires access to all elements of the tensor, and its computational cost scales at least with the size of the tensor, i.e., $\mathcal{O}(N^d)$. In this work, we present \texttt{DEIM-FS} --- a cross algorithm that builds a rank-$r$ Tucker tensor model by sampling $\mathcal{O}(dr^{d-1})$ strategically selected fibers.   The fibers are selected using the DEIM algorithm. The computational cost of \texttt{DEIM-FS} scales linearly with $N$. Our numerical results demonstrate that \texttt{DEIM-FS} remains well-conditioned as $r$ increases and the error of the algorithm closely follows the same-rank HOSVD errors. 

We augment the \texttt{DEIM-FS} algorithm with two capabilities: (i) the introduction of the \texttt{DEIM-FS (iterative)} algorithm, which, unlike \texttt{DEIM-FS}, does not require access to the left singular vectors of the tensor unfolding. Consequently, \texttt{DEIM-FS (iterative)} can be viewed as a black-box Tucker tensor model algorithm. (ii) Additionally, we introduce the rank-adaptive \texttt{DEIM-FS}, where the Tucker rank is adjusted to meet a specified accuracy threshold.

Finally, we introduce \texttt{DLRA-DEIM-FS} to tackle computational cost challenges encountered when solving DLRA equations in nonlinear TDEs. The computational savings of DLRA diminish in cases where the exact rank of the right-hand side tensor ($\mathcal F(\hat{\mathcal{V}})$) is exceptionally high, even when $\hat{\mathcal{V}}$ is low-rank, such as in nonlinear TDEs. We demonstrate that \texttt{DLRA-DEIM-FS} mitigates these cost issues by constructing a low-rank Tucker tensor approximation of the right-hand side tensor in the TDE. Utilizing \texttt{DEIM-FS}, we substantially reduce the required memory and computational cost involved in solving DLRA evolution equations.

\section*{Acknowledgments}
We thank Professor Peyman Givi for many stimulating discussions that led to a number of improvements in this paper. This work is sponsored by the National Science Foundation (NSF), USA under Grant CBET-2152803 and by the Air Force Office of Scientific Research award no. FA9550-21-1-0247. 

\clearpage
\section*{Appendix 1}
The DEIM pseudocode is presented via Algorithm \ref{alg:DEIM}. This algorithm is adopted from \cite{SCNonlinear}. 
\begin{algorithm}
\SetAlgoLined
\KwIn{$\mathbf{U}_{p}=\left[\begin{array}{llll}\mathbf{u}_{1} & \mathbf{u}_{2} & \cdots & \mathbf{u}_{p}\end{array}\right]$}
\KwOut{$\bm I_{p}$}
$\left[\rho, \bm I_{1}\right]=\max \left|\mathbf{u}_{1}\right|$ \hspace{1.3cm} $\rhd$ choose the first index\;
$\mathbf{P}_{1}=\left[\mathbf{e}_{\bm I_{1}}\right]$ \hspace{2.45cm} $\rhd$ construct first measurement matrix\;
\For{$i=2$ \KwTo $p$}{
$\mathbf{P}_{i}^{T} \mathbf{U}_{i} \mathbf{c}_{i}=\mathbf{P}_{i}^{T} \mathbf{u}_{i+1}$ \hspace{0.46cm}  $\rhd$ calculate $c_{i}$\;
$\mathbf{R}_{i+1}=\mathbf{u}_{i+1}-\mathbf{U}_{i} \mathbf{c}_{i}$ \hspace{0.35cm}  $\rhd$ compute residual\;
$\left[\rho, \bm I_{i}\right]=\max \left|\mathbf{R}_{i+1}\right|$ \hspace{0.38cm} $\rhd$ find index of maximum residual\;
$\mathbf{P}_{i+1}=\left[\begin{array}{ll}\mathbf{P}_{i} & \mathbf{e}_{\bm I_{i}}\end{array}\right]$ \hspace{0.38cm} $\rhd$ add new column to measurement matrix\;}
\caption{\texttt{DEIM} Algorithm \cite{SCNonlinear}}
\label{alg:DEIM}
\end{algorithm}

\bibliographystyle{ieeetr}
\bibliography{cas-refs,Hessam}

\end{document}